\documentclass[12pt]{article}
\usepackage{amscd,amssymb,amsmath,latexsym,enumerate}
\usepackage[mathscr]{euscript}
\usepackage{epsfig}
\usepackage{epsfig,psfrag}
\usepackage{here}
\usepackage{graphicx} 
\title{A Spectral Sequence for the $K$-theory of Tiling Spaces\footnote{Work supported by the NSF grants no. DMS-0300398 and no. DMS-0600956.}} 

\author{Jean Savinien, Jean Bellissard\\
{\small Georgia Institute of Technology, School of Mathematics, 
Atlanta GA 30332-0160}
}

\date{ }

\newtheorem{theo}{Theorem}
\newtheorem{defini}{Definition}
\newtheorem{proposi}{Proposition}
\newtheorem{lemma}{Lemma}
\newtheorem{coro}{Corollary}
\newtheorem{rem}{Remark}

\newcommand{\Aa}{{\mathcal A}}
\newcommand{\Bb}{{\mathcal B}}
\newcommand{\Cc}{{\mathcal C}}

\newcommand{\Ff}{{\mathcal F}}

\newcommand{\Kk}{{\mathcal K}}
\newcommand{\Ll}{{\mathcal L}}

\newcommand{\Pp}{{\mathcal P}}

\newcommand{\Ss}{{\mathcal S}}
\newcommand{\Tt}{{\mathcal T}}

\newcommand{\Vv}{{\mathcal V}}

\newcommand{\id}{{\mathbf 1}}

\newcommand{\NM}{{\mathbb N}}
\newcommand{\QM}{{\mathbb Q}}
\newcommand{\PM}{{\mathbb P}}
\newcommand{\RM}{{\mathbb R}}

\newcommand{\SM}{{\mathbb S}}
\newcommand{\TM}{{\mathbb T}}

\newcommand{\ZM}{{\mathbb Z}}

\newcommand{\PG}{{\mathfrak P}}
\newcommand{\pG}{{\mathfrak p}}

\newcommand{\TV}{\Tt_{\PM}}                        
\newcommand{\tra}{\mbox{\sc t}}                    
\newcommand{\Cs}{$C^{\ast}$-algebra }              
\newcommand{\CS}{$C^{\ast}$-algebra}               
\newcommand{\Css}{$C^{\ast}$-algebras }            
\newcommand{\CsS}{$C^{\ast}$-algebras}             
\newcommand{\Ima}{\mbox{\rm Im}}                   
\newcommand{\Ker}{\mbox{\rm Ker}}                  

\newcommand{\qed}{\hfill $\Box$}

\begin{document}

\maketitle

\begin{abstract}
\noindent 
Let $\Tt$ be an aperiodic and repetitive tiling of $\RM^d$ with finite local complexity. 
We present a spectral sequence that converges to the $K$-theory of $\Tt$ with page-$2$ given by a new
cohomology that will be called PV in reference to the Pimsner-Voiculescu exact sequence. It is a generalization of the Serre spectral sequence. The PV cohomology of $\Tt$ generalizes the cohomology of the base space of a fibration with local coefficients in the $K$-theory of its fiber. We prove that it
is isomorphic to the \v{C}ech cohomology of the hull of $\Tt$ (a compactification of the family of its translates).
\end{abstract}




\section{Main Results}
\label{Ktiling07.sec-intro}

\noindent Let $\Tt$ be an aperiodic and repetitive tiling of $\RM^d$ with finite local complexity (definition \ref{Ktiling07.def-repFLC}).
The hull $\Omega$ is a compactification, with respect to an appropriate topology, of the family of translates of $\Tt$ by vectors of $\RM^d$
(definition \ref{Ktiling07.def-hull}). The tiles of $\Tt$ are given compatible $\Delta$-complex decompositions (section \ref{Ktiling07.ssect-ringtrans}), with each simplex punctured,
and the $\Delta$-transversal $\Xi_\Delta$ is the subset of $\Omega$ corresponding to translates of $\Tt$ having the
puncture of one of those simplices at the origin $0_{\RM^d}$.
The prototile space $\Bb_0$ (definition \ref{Ktiling07.def-B0}) is built out of the prototiles of $\Tt$ (translational equivalence classes of tiles)
by gluing them together according to the local configurations of their representatives in the tiling.

\noindent The hull is given a dynamical system structure via the natural action of the group $\RM^d$ on itself by
translation \cite{BBG06}. The \Cs of the hull is isomorphic to the crossed-product \Cs $C(\Omega) \rtimes \RM^d$. 

\noindent There is a map $\pG_o$ from the hull onto the prototile space (proposition \ref{Ktiling07.prop-projhullB0})

\[
\begin{array}{ccl}
\Xi_\Delta & \hookrightarrow & \Omega \\
 & & \downarrow  \pG_0 \\
& & \Bb_0
\end{array}
\]

\noindent which, thanks to a lamination structure on $\Omega$ (remark \ref{Ktiling07.rem-lamination}),
resembles (although is {\it not}) a fibration with base space $\Bb_0$ and fiber $\Xi_\Delta$ (remark \ref{Ktiling07.rem-fibrahull}).
The Pimsner-Voiculescu (PV) cohomology $H_{PV}^\ast$ of the tiling (definition \ref{Ktiling07.def-pimcoho}) is a cohomology of the base space $\Bb_0$
with ``local coefficients'' in the $K$-theory of the fiber $\Xi_\Delta$ (remark \ref{Ktiling07.rem-localcoef}).

\begin{theo}
\label{Ktiling07.thm-SSpim}
There is a spectral sequence that converges to the $K$-theory of the \Cs of the hull
\[
E_2^{rs} \Rightarrow  K_{r+s+d} \bigl( C(\Omega) \rtimes \RM^d \bigr) \,,
\]
\noindent and whose page-$2$ is given by
\[
E_2^{rs} \cong
H^r_{PV} \bigl( \Bb_0 ; K^s (\Xi_\Delta) \bigr) \,. 
\]
\end{theo}

\noindent By an argument using the Thom-Connes isomorphism \cite{Co81} the $K$-theory of $C(\Omega) \rtimes \RM^d$ is isomorphic to the topological $K$-theory
of $\Omega$ (with a shift in dimension by $d$), and the above theorem can thus be seen formally as a generalization of the Serre spectral sequence \cite{Se51} for a
certain class of laminations \( \Xi_\Delta \hookrightarrow \Omega \rightarrow \Bb_0\), which are foliated spaces \cite{MS88} but not fibrations.

\vspace{.1cm}

\noindent This result brings a different point of view on a problem solved earlier by Hunton and Forrest in \cite{FH99}. They built a spectral sequence for the $K$-theory of a crossed product \Cs of a $\ZM^d$-action on a Cantor set. Such an action exists always for tilings of finite local complexity, but it is by no means canonical.
Indeed, thanks to a result of Sadun and Williams \cite{SW03}, the hull of a repetitive tiling with finite local complexity is homeomorphic to a fiber bundle over a torus with fiber
the Cantor set. These two results are sufficient to get the $K$-theory of the hull which, thanks to the Thom-Connes theorem \cite{Co81}, gives also the $K$-theory of the \Cs of the tiling $\Cc(\Omega)\rtimes \RM^d$. 
However, the construction of the hull through an inverse limit of branched manifolds,
initiated by Anderson and Putnam \cite{AP98} for the case of substitution tilings and generalized in \cite{BBG06} to all repetitive tilings with finite local complexity,
suggests a different and more canonical construction. So far however, it is not yet efficient for practical calculations. 

\vspace{.1cm}

\noindent One dimensional repetitive tilings with finite local complexity are all Morita
equivalent to a $\ZM$-action on a Cantor set. The Pimsner-Voiculescu exact sequence \cite{PV80} is then sufficient to compute the $K$-theory of the hull
\cite{Be93}. In the late nineties, before the paper by Forrest and Hunton was written, Mihai Pimsner suggested to one of the
authors\footnote{J.B. is indebted to M. Pimsner for this suggestion} a method to generalize the theorem to $\ZM^d$-actions.
This spectral sequence was used and described already in \cite{BKL01} and is a special case of the Kasparov spectral sequence \cite{Ka88} for $KK$-theory. 
We recall it here for completeness and to justify naming $H_{PV}^\ast$ after Pimsner and Voiculescu.

\begin{theo}
\label{Ktiling07.thm-PimSS2}
\noindent Let $\Aa$ be a \Cs endowed with a $\ZM^d$ action $\alpha$ by $\ast$-automorphisms. The PV complex is defined as
$K_\ast(\Aa)\otimes \Lambda^\ast \ZM^d  \xrightarrow{ \ d_{PV} \ }K_\ast(\Aa)\otimes \Lambda^\ast \ZM^d$ with
\[
d_{PV}= \sum_{i=1}^d (\alpha_{i\ast} - \id)\otimes e_i\wedge \,,
\]
where $\{e_1, \cdots , e_d\}$ is the canonical basis of $\ZM^d$, $\alpha_i= \alpha_{e_i}$ is the restriction of $\alpha$ to the $i$-th component of $\ZM^d$,
whereas $x\wedge$ is the exterior multiplication by $x\in \ZM^d$.

\noindent There is a spectral sequence converging to the $K$-theory of $\Aa$
\[
E_2^{rs} \Rightarrow  K_{r+s+d} \bigl( \Aa \rtimes_\alpha \ZM^d \bigr) \,,
\]
\noindent with page-$2$ isomorphic to the cohomology of the PV complex. 
\end{theo}

\noindent {\it Idea of the proof.} The spectral sequence is built out of the cofiltration associated with the filtration of the mapping torus $M_\alpha(\Aa)$ by ideals of functions vanishing on the squeleton of the torus. The differential for PV cohomology can be identified using the $K$-theory maps (Bott and boundary maps). The reader is referred to \cite{Bla98} and in particular to proposition 10.4.1 where the problem for $d=1$ is treated briefly. \qed

\vspace{.1cm}

\noindent  A more topological expression of this theorem consists in replacing $\ZM^d$ by its classifying space, the torus $\TM^d$,
with a $CW$-complex decomposition given by an oriented (open) $d$-cube and all of its (open) faces in any dimension. Then the PV complex can be proved
to be isomorphic to the following complex: the cochains are given by covariant maps $\varphi (e) \in K_\ast(\Aa)$, where $e$ is a cell of $\TM^d$ and
$\varphi(\overline{e}) =-\varphi(e)$ if $\overline{e}$ is the face $e$ with opposite orientation. The covariance means that if two cells $e,e'$ differ by a
translation $a\in \ZM^d$ then $\varphi(e') = \alpha^a \varphi(e)$. The differential is the usual one, namely $d\varphi (e) = \sum_{e'\in\partial e} \varphi(e')$.
If $H^\ast_P \bigl( \TM^d ; K_\ast (\Aa) \bigr)$ denotes the corresponding cohomology this gives

\begin{coro}
\label{Ktiling07.cor-pim}
The PV cohomology group for the crossed product $\Aa\rtimes_\alpha \ZM^d$ is isomorphic to $H^\ast_{PV} \bigl( \TM^d ; K_\ast (\Aa) \bigr)$.
\end{coro}

\noindent The spectral sequence used by Forrest and Hunton in \cite{FH99} in the case $\Aa=\Cc(X)$ where $X$ is the Cantor set coincides with the PV spectral sequence.
The present paper generalizes this construction for tilings by replacing the classifying space $\TM^d$, by the prototile space $\Bb_0$.

\vspace{.1cm}
\noindent The hull can also be built
out of a {\em box decomposition} \cite{BBG06}. Namely  a {\em box} is a local product of the transversal (which is a Cantor set) by a polyhedron in $\RM^d$ called
the {\em base} of the box. Then the hull can be shown to be given by a finite number of such boxes together with identifications on the boundary of the
bases. This is an extension of the mapping torus, for which there is only one box with base given by a {\em cube} \cite{Se51}. Since $\Aa$ in the present
case is the space of continuous functions on the Cantor set $X$, $K_0(\Aa)$ is isomorphic to $\Cc(X,\ZM)$ whereas $K_1(\Aa)=0$, leading to the present
result. In Section~\ref{Ktiling07.ssect-PimCoho1d} an example of an explicit calculation of the PV cohomology is proposed illustrating the way it
can be used. Further examples, and methods of calculation of PV cohomology will be investigated in future research. However, as it turns out (see section \ref{Ktiling07.ssect-cohom}), the PV cohomology is isomorphic to other cohomologies used so far on the hull, such as the
\v{C}ech cohomology \cite{AP98,Sa03}, the group cohomology \cite{FH99} or the pattern equivariant cohomology \cite{Ke03,KP06,Sa06}.

\vspace{.1cm}

\noindent {\bf Aknowledgements.} Both authors were invited and funded by the Erwin Schr\"odinger Institute (Vienna, Austria) in Summer 2006 and at the Institut Henri Poincar\'e
(Paris, France) in Spring 2007. It is a pleasure for J. Savinien to thank the Erwin Schr\"odinger Institute that supported him as a Junior Research Fellow; most of this work was
done there. He aknowledges financial support from the School of Mathematics of the Georgia Institute of Technology. He thanks J. Hunton and C. Schochet for useful discussions
and comments on this work. J. Bellissard thanks Mihai Pimsner for his explanation of the spectral sequence for $\ZM^d$-actions (theorem \ref{Ktiling07.thm-PimSS2}) provided
in June 1997 in a conversation in Geneva. He is also indebted to J. Hunton, J. Kellendonk and A. Legrand for patiently explaining what a spectral sequence is. He is also
indebted to J.-M. Gambaudo and L. Sadun for explaning him how to built the branched manifolds involved in this work.

\section{Historic Background}
\label{Ktiling07.sec-hist}

\noindent {\em Spectral sequences} were introduced by Leray \cite{Le46} during WWII, as a way to compute the cohomology of a  sheaf. It was later
put in the framework used today by Koszul \cite{Ko47}. One of the first applications of this method was performed by Borel and Serre \cite{BS50}. Later, Serre, in his Ph. D.
Thesis \cite{Se51} defined the notion of Serre fibration and build a spectral sequence to compute the their singular homology or cohomology. This more or less led him to
calculate the rational homotopy of spheres.

\vspace{.1cm}

\noindent In the early fifties, Hirzebruch made an important step in computing the Euler characteristic of various complex algebraic varieties and complex
vector bundles over them \cite{Hi54}. He showed that this characteristic can be computed from the Chern classes of the tangent
bundle and of the vector bundle through universal polynomials \cite{Hi56} which coincides with the Todd genus in the case of varieties. It allowed him to show
that the Euler characteristic is additive for extensions namely if $E,E',E''$ are complex vector bundles and if $0\rightarrow E'\rightarrow E\rightarrow E''\rightarrow 0$
is an exact sequence, than $\chi(E)= \chi(E')+\chi(E'')$. This additivity property led Grothendieck to define axiomatically an additive group characterizing this additivity
relation, which he called the $K$-group \cite{Gro58}. It was soon realized by Atiyah and Hirzebruch \cite{AH61} that the theory could be extended to topological
spaces $X$ and they defined the topological $K$-theory as a cohomology theory without the axiom of dimension. The $K_0$-group is the set of stable equivalence
classes of complex vector bundles over $X$, while $K_1$ (and more generaly $K_n$) is the set of stable equivalence
classes of complex vector bundles over the suspension (respectively the $n$-th suspension)  of $X$. The Bott periodicity theorem reduces the number of groups to
two only. Moreover, the Chern character was shown to define a natural map between the $K$-group and the integer
\v{C}ech cohomology of $X$, and that it becomes an isomorphism when both groups are rationalized. In this seminal paper \cite{AH61} Atiyah and Hirzebruch define a spectral sequence that
will be used in the present work. It is a particular case of Serre spectral sequence for the trivial fibration of $X$ by itself with fiber a point: it converges to the $K$-theory of $X$
and its page-$2$ is isomorphic to the cohomology of $X$.

\vspace{.1cm}

\noindent Almost immediately after this step, Atiyah and Singer extended the work of Hirzebruch to the {\em Index theorem} \cite{AS63,AS68}
for elliptic operators. Such an operator is defined between two vector bundles, it is unbounded, in general and, with the correct domains of definition defines
a Fredholm operator. The index can be interpreted as an element of   the $K$-group and is calculated through a formula which generalizes the results of
Hirzebruch for algebraic varieties. Eventually, Atiyah and Singer extended the theory to the equivariant $K$-theory valid if a compact Lie group $G$ acts on the
vector bundle. If $P$ is an elliptic operator commuting with $G$ its index gives an element of the covariant $K$-theory. In a programmatic paper \cite{Si70},
Singer proposed various extensions to elliptic operators with c{\oe}fficients depending on parameters. As an illustration of such a program, Coburn, Moyer
and Singer \cite{CMS73} gave an index theorem for elliptic operators with almost periodic c{\oe}fficients (see also \cite{Su79}). Eventually the index theorem
became the cornerstone of Connes' program to build a {\em Noncommutative Geometry}. In a seminal paper \cite{Co79} he defined the theory of
{\em noncommutative integration} and showed that its first application was an index theorem for elliptic operators on a foliation (see also \cite{Co82,Co94}). 

\vspace{.1cm}

\noindent In the seventies, it was realized that the Atiyah-Hirzebruch $K$-theory could be expressed in algebraic terms through the \Cs $\Cc(X)$ of continuous functions on the compact space $X$.
The definition of the $K$-group requires then to consider matrix valued continuous functions $M_n(\Cc(X))$ instead
for all $n$. The smallest \Cs containing all
of them is $\Cc(X)\otimes \Kk$, where $\Kk$ denotes the \Cs of compact operators on a separable Hilbert space. Therefore, since this later algebra is non commutative, all the construction could be
used for any \CS. Then Kasparov \cite{Ka80, Ka88} defined the notion of $KK$-theory generalizing even more the $K$ groups to correspondences between two \CsS.

\vspace{.1cm}

\noindent The problem investigated in the present paper is the latest development of a program that was initiated in the early eighties \cite{JM82,Be82}
when the first version of the {\em gap labeling theorem} was proved
(see \cite{Be86,Be93,BBG06} for later developments). At that time the problem was to compute the spectrum of a Schr{\"o}dinger operator $H$ in an aperiodic potential. Several examples where discovered
of Schr\"odinger operators with a Cantor-like spectrum \cite{Hof76,Mos81}. Labeling the infinite number of gaps per unit length interval, was a challenge. It was realized that the $K$-theory class of the
spectral projection on spectrum below the gap was a proper way of doing so \cite{Be82}. The first calculation was made on the Harper equation and gave an explanation for a result already obtained by
Claro and Wannier \cite{CW78}, a result eventually used in the theory of the quantum Hall effect \cite{TKN2}. This problem was motivated by the need for a theory of aperiodic solids, in particular their
electronic and transport properties. With the discovery of {\em quasicrystals} in 1984 \cite{SBGC84}, this question became crucial in Solid State Physics. The construction of the corresponding
\Cs became then the main issue and led to the definition of the {\em Hull}  \cite{Be86}. It was proved that the Hull is a compact metrizable space $\Omega$ endowed with an action of $\RM^d$ via
homeomorphisms. Then it was proved in \cite{Be86,Be93}, that the resolvent of the Schr\"odinger operator $H$ belongs to the \Cs $\Aa =\Cc(\Omega)\rtimes \RM^d$. With each $\RM^d$-invariant
probability measure $\PM$ on $\Omega$ is associated a trace $\TV$ on this algebra. Following an argument described in \cite{Be82}, it was proved that a gap could be labeled by the value of the
density of state, and that this value belongs to the image by the trace $\TV$ of the group $K_0(\Aa)$. During the eighties several results went on to compute the set of gap labels \cite{Pi83,Be86}.
In one dimension, detailed results could be proved (see \cite{Be93} for a review). The most spectacular result was given in the case of a discretized Laplacian with a potential taking on finitely many
values (called {\em letters}): in such a case the set of gap labels is the $\ZM$-module generated by the occurrence probability of all possible finite words found in the sequence defined by the
potential. If this sequence is given by a {\em substitution}, these occurrence probabilities can be computed explicitely in terms of the incidence matrix of the substitution and of the associated
substitution induced on the set of words with two letters \cite{QU87,Be93}. The key property in proving such results was the use of the Pimsner-Voiculescu exact sequence \cite{PV80}. Soon after,
A. van Elst \cite{vEl94} extended these results to the case of 2D-potentials, using the same method.

\vspace{.1cm}

\noindent In the French version \cite{Co90} of his book on {\em Noncommutative geometry} \cite{Co94}, Connes showed how the general formalism he had developed could be illustrated
with the special case of the Penrose tiling. Using the substitution rules for its construction, he introduced a \CS, which turns out to be AF, and computed its ordered $K_0$-group, using the classification
of AF-algebras obtained by Bratelli in 1972 \cite{Br72}. This result was an inspiration for Kellendonk, who realized that, instead of looking at the inflation rule as a source of noncommutativity, it was
actually better to consider the space translations of the tiling \cite{Ke95}, in the spirit of the formalism developed in \cite{Be82,Be86,Be93} for aperiodic solids. He extended this latter construction of
the Hull to tilings and then this Hull is called {\em tiling space}. In this important paper, Kellendonk introduced the notion of {\em forcing the border} for a general tiling, which appears today as an important
property for the calculation of the K-groups and cohomology of a tiling space. This work gave a strong motivation to prove the {\em Gap Labeling Theorem in higher dimension} and to compute the $K$-group
of the Hull. It was clear that the method used in one dimension, through the Pimsner-Voiculescu exact sequence could only be generalized through a spectral sequence. The first use of spectral sequences
in computing the set of gap labels on the case of the 2D-octagonal tiling \cite{BCL98} was followed by a proof of the Gap Labeling Theorem for 3D-quasicrystals \cite{BKL01}. Finally, the work by Forrest and
Hunton \cite{FH99} used a classical spectral sequence to compute the full $K$-theory of the \Cs for an action of $\ZM^d$ on the Cantor set. It made possible the computation of the K-theory and the
cohomology of the Hull for quasicrystals in two and three dimensions \cite{HFK02,GHK05}. Other examples of tilings followed later \cite{Sa05}.

\vspace{.1cm}

\noindent In addition to Kellendonk's work, several important contributions helped to build tools to prove the higher dimensional version of the Gap
Labeling Theorem. Among the major contributions was the work of Lagarias \cite{La99A,La99B,LP03}, who introduced a geometric and combinatoric
aspect of tiling through the notion of a Delone set. This concept was shown to be conceptually crucial in describing aperiodic solids \cite{BHZ00}. Through
the construction of Voronoi, Delone sets and tilings become equivalent concepts, allowing for various intuitive point of views to study such problems.

\vspace{.1cm}

\noindent Another important step was performed in 1998 by Anderson and Putnam \cite{AP98}, who proposed to build the tiling space of a substitution tiling
through a CW-complex built from the prototiles of a tiling (called prototile space in this paper). The substitution induces a map from this CW-complex into itself and they showed that the inverse
limit of such system becomes homeomorphic to the tiling space. A similar construction was proposed independently by Gambaudo and Martens in 1999 to
describe dynamical system. This latter case corresponds to 1D repetitive tilings with finite local complexity, so that this latter construction goes beyond the
substitution tilings. One interesting outcome of this work was a systematic construction of minimal dynamical systems, uniquely ergodic or not, and with positive
entropy\footnote{This result circulated as a preprint in 1999 but was published only in 2006} (see \cite{GM06}). Eventually, the Gambaudo-Martens construction
led to the construction of the Hull as an inverse limit of compact oriented branched flat Riemannian manifolds
\cite{BBG06}\footnote{This paper was posted on {\tt arXiv.com math.DS/0109062} in its earlier version in 2001 but was eventually publishhed in a final form in
2006 only.} used in the present paper (although only their topological $CW$-complex structure is needed here).
Equivalently, the Hull can be seen as a {\em lamination} \cite{Gh99} of a {\em foliated space} \cite{MS88}. The extension of this construction to include
tilings without finite local complexity, such as the pinwheel model \cite{Rad99}, was performed by Gambaudo and Benedetti \cite{BG03} using the notion of solenoids
(reintroduced by Williams in the seventies \cite{Wil74}).

\section{A Mathematical Reminder}
\label{Ktiling07.sect-mathback}

\noindent The preliminary definitions and results required for stating theorem \ref{Ktiling07.thm-SSpim}
are presented here. They are taken mostly from previous work of the first author
in \cite{BHZ00,BBG06} on tilings and Delone sets and the reader is referred to those papers
for complete proofs.
Let $\RM^d$ denote the usual Euclidean space of dimension $d$ with euclidean norm \( \| \cdot \| \).
First of all the definition of a repetive tiling of $\RM^d$ with finite local complexity is recalled, and then the hull and
transversal of such a tiling are introduced. The connection with Delone sets is briefly mentioned as well as the
definition of the groupoid of the transversal of a Delone set.
The prototile space of a tiling and the ring of functions on the transversal are then built, and finally the
PV cohomology defined.

\subsection{Tilings and their hulls}
\label{Ktiling07.ssect-tiling}

\begin{defini}
\label{Ktiling07.def-tiling}

\begin{enumerate}[(i)]

\item A {\rm tile} of $\RM^d$ is a compact subset of $\RM^d$ which is homeomorphic to the unit ball.

\item A punctured tile is an ordered pair
consisting of a tile and one of its points.

\item A {\rm tiling} of $\RM^d$ is a covering of $\RM^d$ by a family of tiles whose interiors are pairwise disjoint.
A tiling is said to be punctured if its tiles are punctured.

\item A {\rm prototile} of a tiling is a translational equivalence class of tiles (including the puncture).

\item The {\rm first corona} of a tile in a tiling $\Tt$ is the union of the tiles of $\Tt$ intersecting it.

\item A {\rm collared prototile} of $\Tt$ is the subclass of a prototile whose representatives have the same first corona up to translation.

\end{enumerate}

\end{defini}

\noindent A collared prototile is a prototile where a local configuration of its representatives has been specified: each
representative has the same neighboring tiles.

\noindent In the sequel it is implicitely assumed that tiles and tilings are punctured. All tiles are assumed to be finite
$\Delta$-complexes which are particular $CW$-complex structures (see section \ref{Ktiling07.ssect-ringtrans}). They are also required to be compatible
with the tiles of their first coronas, {\it i.e.} the intersection of any two tiles is itself a sub-$\Delta$-complex of both tiles.
In other words, tilings considered here are assumed to be $\Delta$-complexes of $\RM^d$.

\begin{defini}
\label{Ktiling07.def-patch}
Let $\Tt$ be a tiling of $\RM^d$.

\begin{enumerate}[(i)]
\item A {\rm patch} of $\Tt$ is a finite union of neighboring tiles in $\Tt$ that is
homeomorphic to a ball. A patch is punctured by the puncture of one of the tiles that it contains. The radius of a patch is the radius of the smallest ball that contains it.
\item A {\rm pattern} of $\Tt$ is a translational equivalence class of patches of $\Tt$.
\item The {\rm first corona} of a patch of $\Tt$ is the union of the tiles of $\Tt$ intersecting it.
\item A {\rm collared pattern} of $\Tt$ is the subclass of a pattern whose representatives have the same first corona up to translation.
\end{enumerate}

\end{defini}

\noindent The following notation will be used in the sequel: prototiles and patterns will be written with a hat, for instance $\hat{t}$ or
$\hat{p}$, to distinguish them from their representatives. Often the following convention will be implicit: if $\hat{t}$ is a prototile and $\hat{p}$ a
pattern, then $t$ and $p$ will denote their respective representatives that have their punctures at the origin $0_{\RM^d}$.

\noindent The results in this paper are valid for the class of tilings that are {\it repetitive} with {\it finite local complexity}.

\begin{defini} 
\label{Ktiling07.def-repFLC}
Let $\Tt$ be a tiling of $\RM^d$.

\begin{enumerate}[(i)]

\item $\Tt$ has {\rm Finite Local Complexity} (FLC) if for any $R>0$ the set of patterns of $\Tt$ whose representatives have radius less than $R$ is finite.

\item $\Tt$ is {\rm repetitive} if for any patch of $\Tt$ and every $\epsilon>0$, there is an $R>0$ such that for every $x$ in $\RM^d$
there exists {\rm modulo an error} $\epsilon$ with respect to the Hausdorff distance, a translated copy of this patch belonging to $\Tt$ and contained in the ball $B(x,R)$.

\item For any $x$ in $\RM^d$, let \( \Tt +x = \{ t+x \, : \, t \in \Tt \}\) denote its translation, then
$\Tt$ is {\rm aperiodic} if there is no $x\ne0$ in $\RM^d$ such that \( \Tt+x = \Tt\).

\end{enumerate}

\end{defini}

\noindent For tilings with FLC the repetitivity condition (ii) above can be stated more precisely: {\it A tiling $\Tt$ with FLC is repetitive if given any patch there is an
$R>0$ such that for every $x$ in $\RM^d$ there exists an {\em exact copy} of this patch in $\Tt$ contained in the ball $B(x,R)$.}

\vspace{.1cm}

\noindent The class of repetitive tilings satisfying the FLC property is very rich and has been investigated for decades. It started in the 70's
with the work of Penrose \cite{Pen74} and Meyer \cite{Mey72} and went on 
both from an abstract mathematical level and with a view towards applications, in particular to the physics of quasicrystals \cite{La99A,La99B}.
It contains an important subclass of the class of {\it substitution tilings} that was reinvestigated in the 90's by Anderson and Putnam in \cite{AP98},
and also the class of tilings obtained by the {\it cut-and-projection} method, for which a comprehensive study by Hunton, Kellendonk and Forrest can
be found in \cite{HFK02}, and more generaly it contains the whole class of {\it quasiperiodic tilings}, which are models for quasicrystals.

\noindent The Ammann aperiodic tilings are examples of substitution tilings (see \cite{GruShep87} chapter 10, Ammann's original work does not appear in the literature).
In addition other examples include the octogonal tiling as well as the famous Penrose ``kite and darts'' \cite{Pen74} tilings which are both substitution and cut-and-projection
tilings  (see \cite{GruShep87} chapter 10, and \cite{Rad99} chapter 4).
However the so-called Pinwheel tiling (see \cite{Rad94} and \cite{Rad99} chapter 4) is a substitution tiling but
does not satisfy the FLC property given here since prototiles are defined here as equivalence classes of
tiles under only translations and not more general isometries of $\RM^d$ like rotations.

\vspace{.1cm}

\noindent The first author has proposed in \cite{BHZ00} a topology that applies to a large class of tilings (for which there exist an $r_0>0$ such that all tiles contain a ball of
radius $r_0$).
In the case considered here, where tiles are assumed to be finite $CW$-complexes, this topology can be adapted as follows.

\noindent Let $\Ff$ be a family of tilings whose tiles contain a ball of a fixed radius $r_0>0$ and have compatible $CW$-complex structures (the intersection of any two tiles is
a subcomplex of both).
Given an open set $O$ in $\RM^d$ with compact closure and an $\epsilon>0$, a neighborhood of a tiling $T$ in $\Ff$ is given by

\[
U_{O,\epsilon}(T) =
\left\{
T' \in \Ff \, : \,
\sup_{0\le k \le d} h_k ( O \cap T^k, O \cap T'^k) < \epsilon
\right\} \,,
\]
where $T^k$ and $T'^k$ are the $k$-skeletons of $T$ and $T'$ respectively and $h_k$ is the $k$-dimensional Hausdorff distance.

\noindent Let $\Tt$ be a tiling of $\RM^d$. The group $\RM^d$ acts on the set of all translates of $\Tt$, the action (translation)
is denoted \( \tra^a, a \in \RM^d\): \(\tra^a \Tt =\Tt+a =  \{ t+a \, : t \in \Tt \} \).

\begin{defini}
\label{Ktiling07.def-hull}

\begin{enumerate}[(i)]

\item The {\em hull} of $\Tt$, denoted $\Omega$, is the closure of  \( \tra^{\RM^d} \Tt\).

\item The canonical {\em transversal}, denoted $\Xi$, is the subset of $\Omega$ consisting of tilings that have the puncture
of one of their tiles at the origin $0_{\RM^d}$.

\end{enumerate}

\end{defini}

\noindent The hull of a tiling is seen as a dynamical system \( (\Omega, \RM^d, \tra)\) which, for the class of tilings considered here,
has interesting properties that are now stated (see \cite{BBG06} section 2.3).  

\begin{theo} {\rm  \cite{BBG06,La99A}}
\label{Ktiling07.thm-hull}
Let $\Tt$ be a tiling of $\RM^d$.

\begin{enumerate}[(i)]
\item $\Tt$ is repetitive if and only if the dynamical system of its hull \( (\Omega, \RM^d, \tra) \) is {\rm minimal}.
\item If $\Tt$ has FLC, then its canonical transversal $\Xi$ is totaly disconnected.
\item If $\Tt$ is aperiodic, repetitive and has FLC, then $\Xi$ is a Cantor set (perfect and totaly disconnected).
\end{enumerate}

\end{theo}

\noindent The minimality of the hull allows one to see any of its points as just a translate of $\Tt$.

\begin{rem}
\label{Ktiling07.rem-topohull}
{\em A metric topology for tiling spaces has been used in the literature for historical reasons. Let $\Tt$ be a repetive tiling or $\RM^d$ with FLC.
The orbit space of $\Tt$ under translation by vectors of
$\RM^d$, \( \tra^{\RM^d} \Tt\), is endowed with a metric as follows (see \cite{BBG06} section 2.3).
For $T$ and $T'$ in \( \tra^{\RM^d} \Tt \), let $A$ denote the set of $\varepsilon$ in $(0,1)$ such that there exist
$x$ and $x'$ in \( B(0, \varepsilon ) \) for which \( \tra^x T \) and \( \tra^{x'} T'\) agree on \(  B(0, \frac{1}{\varepsilon} ) \),
{\it i.e.} their tiles whose punctures lie in the ball are matching, then

\[
\delta (T,T') = \left\{
\begin{array}{lcl}
\inf A & {\rm if} & A \ne \emptyset \,, \\
1 & {\rm if} &  A = \emptyset \,.
\end{array}
\right.
\]

\noindent Hence the diameter of  \( \tra^{\RM^d} \Tt\) is bounded by $1$ and the action of $\RM^d$ is continuous. 

\noindent For the class of repetitive tilings with FLC, the topology of the hull given in definition \ref{Ktiling07.def-hull} is equivalent to this $\delta$-metric topology \cite{BBG06}.}
\end{rem}

\subsection{Delone sets and groupoid of the transversal}
\label{Ktiling07.ssect-delone}

\noindent The notions for tilings given in the previous subsection can be translated for the set of their punctures in 
terms of Delone sets.

\begin{defini}
\label{Ktiling07.def-delone}
Let $\Ll$ be a discrete subset of $\RM^d$.

\begin{enumerate}[(i)]
\item  Given $r>0$, $\Ll$ is {\rm $r$--uniformly discrete} if any 
open ball of radius $r$ in $\RM^d$ meets $\Ll$
in at most one point.

\item Given $R>0$, $\Ll$ is {\rm $R$--relatively dense} if any open 
ball of radius $R$ in $\RM^d$ meets $\Ll$
in at least one point.

\item $\Ll$  is an {\rm $(r,R)$--Delone set} is it is $r$--uniformly discrete 
and $R$--relatively dense.

\item $\Ll$ is repetitive if given any 
{\em finite} subset $p \subset \Ll$ {\em and any $\epsilon >0$}, there is an 
$R>0$ such that the intersection of $\Ll$ with any closed ball of radius $R$ 
contains a copy (translation) of $p$ {\em modulo an error of $\epsilon$ 
(w.r.t. the Hausdorff distance)}

\item A patch of radius $R>0$ of $\Ll$ is a subset of $\RM^d$ of the form 
\( \left( \Ll - x \right) \cap B(0,R) \), for some $x \in \Ll$.
If for all $R>0$ the set of its patches 
of radius $R$ is finite then $\Ll$ has finite local complexity.

\item $\Ll$ is {\rm aperiodic} if there is no $x\ne0$ in $\RM^d$ such that $\Ll - x  = \Ll$.

\end{enumerate}
\end{defini}

\noindent Condition {\it (iv)} above is equivalent to saying that $\Ll - \Ll$, the set of vectors of $\Ll$, is discrete.

\vspace{.1cm}

\noindent Given a punctured tiling $\Tt$ if there exists $r, R>0$ such that each of its tiles contains a ball a radius $r$ and is
contained in a ball of radius $R$, then its set of punctures $\Ll_\Tt$ is an $(r,R)$-Delone set, and it is repetitive and has FLC
if and only if the tiling is repetitive and has FLC. Conversely, given a Delone set the Voronoi construction below gives a tiling,
and they both share the same repetitivity or FLC properties.

\begin{defini}
\label{Ktiling07.def-voronoi}
Let $\Ll$ be an $(r,R)$--Delone set of $\RM^d$. The Voronoi tile at \( x \in \Ll\), is defined by
\[
T_x = \left\{ y \in \RM^d \, : \, \|y-x\| \le \| y-x'\|, \, \forall x' \in \Ll \right\} \,,
\]
with puncture the point $x$.
The Voronoi tiling $\Vv$ associated with $\Ll$ is the tiling of $\RM^d$ whose 
tiles are the Voronoi tiles of $\Ll$.
\end{defini}

\noindent The tiles of the Voronoi tiling of a Delone set are (closed) convex polytopes that touch on common faces.

\vspace{.1cm}

\noindent The hull and transversal of a Delone set are defined as follows \cite{BHZ00}.

\begin{defini}
\label{Ktiling07.def-delonehull}
Let $\Ll$ be a Delone set of $\RM^d$.
The hull $\Omega$ of $\Ll$ is the closure in the weak--$\ast$ 
topology of the set of translations of the Radon
measure that has a Dirac mass at each point of $\Ll$
\[
\Omega = \overline{ \bigl\{ \tra^a \mu\,;\, a\in \RM^d \bigr\} }^{\, w^\ast} \,,
\quad \quad \mu = 
\sum_{x\in \Ll} \delta_x \,, 
\]
and the transversal $\Xi$ of  $\Ll$ is the subset 
\[ 
\Xi = \bigl\{ \omega \in \Omega \, ; \,     \omega(\{0\}) =1 \bigr\} \,.
\]
\end{defini}

\noindent Given a general Delone set this topology of its hull is strictly coarser than the $\delta$-metric topology for its Voronoi tiling given in remark \ref{Ktiling07.rem-topohull},
but if it is repetitive and has FLC then they are equivalent \cite{BBG06}. The weak-$\ast$ topology is here also
equivalent to the local Hausdorff topology on the set $\tra^{\RM^d} \Ll$ of translates of $\Ll$: given an open set $O$ in $\RM^d$ with compact closure
and an $\epsilon>0$, a neighborhood of  $\ell$ is given by

\[
U_{O,\epsilon}(\ell) = 
\bigl\{
\ell' \in \tra^{\RM^d} \Ll \, : \, h_d ( \ell \cap O , \ell' \cap O ) < \epsilon
\bigr\} \,,
\]

\noindent where $h_d$ is the $d$-dimensional Hausdorff distance.

\vspace{.1cm}

\noindent The hull of a Delone set is also seen as a dynamical system under the homeomorphic action of
$\RM^d$ by translation. Just as in proposition \ref{Ktiling07.thm-hull}, the hull and transversal of repetitive
Delone sets with FLC have similar interesting properties.

\begin{theo} {\rm  \cite{BHZ00}}
\label{Ktiling07.thm-delonehull}
Let $\Ll$ be a Delone set of $\RM^d$.
\begin{enumerate}[(i)]
\item $\Ll$ is repetitive if and only if the dynamical system of its hull \( (\Omega, \RM^d, \tra) \) is {\rm minimal}.
\item If $\Ll$ has FLC, then its transversal $\Xi$ is totally disconnected.
\item If $\Ll$ is aperiodic, repetitive and has FLC, then $\Xi$ is a Cantor set (perfect and totaly disconected).
\end{enumerate}
\end{theo}

\noindent The two approches of tilings or Delone sets which are repetitive and have FLC are thus equivalent. Such tilings give
rise to Delone sets (their sets of punctures), and conversely such Delone sets give rise to
tilings (their Voronoi tilings) and their respective hulls share the same properties.

\vspace{.1cm}

\noindent Recall that a groupoid is a small category in which every morphism is invertible \cite{Co79, Re80}.

\begin{defini} 
\label{Ktiling07.def-groupoid}
Let $\Ll$ be a Delone set of $\RM^d$. The {\em groupoid of the transversal} is the groupoid $\Gamma$ whose
set of objects is the transversal: \(\Gamma^0 = \Xi\), and whose set of arrows is
\[
\Gamma^1 =
\bigl\{
(\xi,x) \in \Xi \times \RM^d \, : \, \tra^{-x}\xi \in \Xi 
\bigr\} \,.
\]

\noindent Given an arrow $\gamma = (\xi,x)$ in $\Gamma^1$, its {\em source} is the object $s(\gamma)= \tra^{-x} \xi \in \Xi$
and its {\em range} the object $r(\gamma) = \xi \in \Xi$.
\end{defini}

\noindent The Delone sets considered here are repetitive with FLC and the groupoids of their transversals are {\it \'etale},
{\it i.e.} given any arrow $\gamma$ with range (or source) object $x$, there exists an open neighborhood $O_x$ of $x$ in $\Gamma^0$ and a homeomorphism
\( \varphi : O_x \rightarrow \Gamma^1\) that maps $x$ to $\gamma$.
This implies that the sets of arrows having any given object as a source or range are discrete. This follows because for such Delone sets
their sets of vectors are discrete and given any fixed $l>0$ their sets of vectors of length less than $l$ are finite.

\subsection{Prototile space of a tiling}
\label{Ktiling07.ssect-boxhull}

\noindent Let $\Tt$ be an aperiodic and repetitive tiling of $\RM^d$ with FLC and assume its tiles are compatible finite $CW$-complexes
(the intersection of two tiles is a subcomplex of both). A finite $CW$-complex $\Bb_0$, called prototile space, is built out of its prototiles by
gluing them along their boundaries according to all the local configurations of their representatives in $\Tt$.

\begin{defini}
\label{Ktiling07.def-B0}
Let $\hat{t}_j , \, j=1, \cdots N_0$, be the prototiles of $\Tt$.  Let $t_j$ denote the representative of $\hat{t}_j$ that has its puncture at the origin.
The {\em prototile space} of $\Tt$, \(\Bb_0(\Tt)\),  is the quotient $CW$-complex
\[ 
\Bb_0 (\Tt) =  \coprod_{j=1}^{N_0} t_j \; / \sim \,,
\] 

\noindent where two $n$-cells $e_i^n \in t^n_i$ and $e_j^n \in t^n_j$ are identified if there exists $u_i, u_j \in \RM^d$ for which
\(t_i + u_i\) and \(t_{j}+u_{j}\) are tiles of $\Tt$ such that \(e_i^n + u_i\) and \(e_j^n + u_j\) coincide on the intersection of
their $n$-skeletons.

\noindent The {\em collared prototile space} of $\Tt$, \(\Bb_0^c (\Tt)\), is built similarly from the collared prototiles of $\Tt$: \(t_i^c, i=1, \cdots N_0^c\).

\end{defini}

\noindent The images in \(\Bb_0(\Tt)\) or \(\Bb_0^c (\Tt)\) of the tiles $t_j^{(c)}$'s will be denoted $\tau_j$ and still be called tiles.

\begin{proposi} 
\label{Ktiling07.prop-projhullB0}
There is a continuous map \( \pG_{0,\Tt}^{(c)} : \Omega \rightarrow \Bb_0^{(c)}(\Tt) \) from the hull onto the (collared) prototile space.
\end{proposi}

\noindent {\it Proof.} Let \( \lambda_0^{(c)} :  \coprod_{j=1}^{N_0^{(c)}} t_j \rightarrow {\mathcal B}_0^{(c)}(\Tt) \) be the quotient map.
And let \( \rho_0^{(c)} : \Omega \times \RM^d \rightarrow \coprod_{j=1}^{N_0^{(c)}} t_j^{(c)} \) be defined as follows.
If $x$ belongs to the intersection of $k$ tiles \(t^{\alpha_1}, \cdots t^{\alpha_k}\),
in $\omega$, with \(t^{\alpha_l} = t_{j_l}^{(c)} + u_{\alpha_l}(\omega), \, l=1, \cdots k\), then
\( \rho_0^{(c)} (\omega,x) =  \coprod_{l=1}^{k} x-u_{\alpha_l}(\omega) \) and lies in the disjoint union of the $t_{j_l}^{(c)}$'s.

\noindent The map $\pG_{0,\Tt}^{(c)}$ is defined as the composition: \( \omega  \mapsto  \lambda_0^{(c)} \circ \rho_0^{(c)} ( \omega, 0_{\RM^d} ) \).
The map \( \rho_0^{(c)} (\cdot, 0_{\RM^d})\) sends the origin of $\RM^d$, that lies in some tiles of $\omega$,
to the corresponding tiles $t_j^{(c)}$'s at the corresponding positions.

\vspace{.1cm}

\noindent In \(\Bb_0^{(c)}(\Tt)\), points on the boundaries of two tiles $\tau_j^{(c)}$ and $\tau_{j'}^{(c)}$ are identified if there are neighboring copies
of the tiles $t_j^{(c)}$ and $t_{j'}^{(c)}$ somewhere in $\Tt$ such that the two associated points match.
This ensures that the map $\pG_{0,\Tt}^{(c)}$ is well defined, for if in $\RM^d$ tiled by $\omega$, the origin
$0_{\RM^d}$ belongs to the boundaries of some tiles, then the corresponding points in \(\coprod_{j=1}^{N_0^{(c)}} t_j^{(c)} \)
given by \(\rho_0^{(c)} (\omega, 0_{\RM^d})\) are identified by $\lambda_0^{(c)}$.

\vspace{.1cm}

\noindent Let $x$ be a point in \( \Bb_0^{(c)}(\Tt)\), and $O_x$ and open neighborhood of $x$. Say $x$ belongs to the intersection of some tiles
\(\tau_{j_1}^{(c)}, \cdots \tau_{j_k}^{(c)}\). Let $\omega$ be a preimage of $x$: \(\pG_{0,\Tt}^{(c)} (\omega) = x\).
The preimage of $O_x$ is the set of tilings $\omega'$'s for which the origin lies in some neighborhood of tiles that are translates of
\( t_{j_1}^{(c)}, \cdots t_{j_k}^{(c)}\), and this defines a neighborhood of $\omega$ in the hull. Therefore \(\pG_{0,\Tt}^{(c)}\) is continuous. 
\qed

\vspace{.1cm}

\noindent For simplicity, the prototile space \(\Bb_0(\Tt)\) is written \(\Bb_0\), and the map \(\pG_{0,\Tt}\) is written \(\pG_0\).

\vspace{.1cm}

\noindent The lift of the puncture of the tile $\tau_j$ in $\Bb_{0}$, denoted $\Xi(\tau_j)$, is a subset of the transversal called
the {\it acceptance zone} of the prototile $\hat{t}_j$. It consists of all the tilings that have the puncture of a representative of $\hat{t}_j$ at the
origin. The $\Xi(\tau_j)$'s for $j=1, \cdots N_0$, form a clopen partition of the transversal, because any element of $\Xi$ has the puncture of a unique tile at the origin which
corresponds to a unique prototile. The $\Xi(\tau_j)$'s  are thus Cantor sets like $\Xi$.

\begin{rem}
\label{Ktiling07.rem-lamination}
{\em
Although those results will not be used here, it has been proven in \cite{BBG06} that the prototile space $\Bb_0$ (as well as the
patch spaces $\Bb_p$ defined similarly in the next section, definition \ref{Ktiling07.def-Bp}) has the structure of a flat oriented Riemanian branched manifold.
Also, the hull $\Omega$ can be given a {\em lamination} structure as follows: the lifts of the interiors of the tiles
\( B_{0j} = \pG_0^{-1}( \stackrel{\circ}{\tau_j} )\) are boxes
of the lamination which are homeomorphic to \( \stackrel{\circ}{t_j} \times \Xi( \tau_j )\) via the maps \( (x,\xi) \mapsto \tra^{-x}\xi \) which
read as local charts.}
\end{rem}

\subsection{The hull as an inverse limit of patch spaces}
\label{Ktiling07.ssect-invlim}

\noindent As in the previous section, let $\Tt$ be an aperiodic and repetitive tiling of $\RM^d$ with FLC, and assume its tiles are compatible
finite $CW$-complexes (the intersection of two tiles is a subcomplex of both). 
Let $\Ll_\Tt$ denote the Delone set of punctures of $\Tt$. It is repetitive and has FLC.
Let  $\Pp_\Tt$ denote the set of patterns of $\Tt$. As $\Tt$ has finite local complexity (hence finitely many prototiles),
the set $\Pp_\Tt$ is countable, and for any given $l>0$ the set of patterns whose representatives have radius less than $l$ is finite.

\vspace{.1cm}

\noindent A finite $CW$-complex $\Bb_p$, called a patch space, associated with a pattern $\hat{p}$ in $\Pp_\Tt$ is built from the prototiles of an
appropriate subtiling of $\Tt$,  written $\Tt_p$ below, in the same way that $\Bb_0$ was built from the prototiles of $\Tt$ in definition \ref{Ktiling07.def-B0}. 
The construction goes as follows. 

\vspace{.1cm}

\noindent Let $\hat{p}$ in  $\Pp_\Tt$ be a pattern of $\Tt$.

\begin{enumerate}[(i)]

\item Consider the sub-Delone set $\Ll_p$ of $\Ll_\Tt$ consisting of punctures of all the representative patches in $\Tt$ of $\hat{p}$. 
$\Ll_p$ is repetitive and has FLC. 

\item The Voronoi tiling  $\Vv_p$ of $\Ll_p$ is built and each point of $\Ll_\Tt$ is assigned to a unique tile of $\Vv_p$ as explained below.

\item Each tile $v$ of $\Vv_p$ is replaced by the patch $p_v$ of $\Tt$ made up of the tiles whose punctures
have been assigned to $v$. This gives a repetitive tiling with FLC, $\Tt_p$, whose tiles are those patches $p_v$'s.

\item $\Bb_p$ is built out of the collared prototiles of $\Tt_p$, by gluing them along their boundaries 
according to the local configurations of their representatives in $\Tt_p$.

\end{enumerate}

\noindent  The second point needs clarifications since the tiles of $\Vv_p$ are Voronoi tiles  (convex polytopes, see definition \ref{Ktiling07.def-voronoi})
of $\Ll_p$ and not patches of $\Tt$.  If a point of $\Ll_\Tt$ (a puncture of a tile of $\Tt$) lies on the boundary of some (Voronoi) tiles of $\Vv_p$, a
criterion for assigning it to a specific one is required.
To do so, let $u$ be a vector of $\RM^d$ that is not colinear to any of the faces of the tiles of $\Vv_p$ (such a vector
exists since $\Vv_p$ has FLC, hence finitely many prototiles). A point $x$ is said to be $u$-interior to a subset $X$ of $\RM^d$ if there exist an
$\epsilon>0$ such that \(x+\epsilon u\) belongs to the interior of $X$. Since $u$ is not colinear to any of the faces of the tiles of $\Vv_p$, if a point
$x$ belongs to the intersection of the boundaries of several (Voronoi) tiles of $\Vv_p$, it is $u$-interior to only one of them.
This allows as claimed in (ii) to assign each point of $\Ll_\Tt$ to a unique tile of $\Vv_\Tt$.
Now as explained in (iii), each Voronoi tile $v$ can then be replaced by the patch of $\Tt$ which is the union of
the tiles of $\Tt$ whose punctures are $u$-interior to $v$.

\vspace{.1cm}

\noindent Each patch $p_v$ is considered a tile of $\Tt_p$ and punctured by the puncture of the Voronoi tile
$v$ which is by construction the puncture of some representative patch of $\hat{p}$ in $\Tt$. As patches of $\Tt$, the $p_v$'s are also
compatible finite $CW$-complexes as they are made up of tiles of $\Tt$ which are.

\vspace{.1cm}

\noindent The prototiles of $\Tt_p$ are actually patterns of $\Tt$, and thus $\Tt_p$ is considered
a subtiling of $\Tt$. From this remark it can be proven that there is a homeomorphism between the hull of $\Tt_p$
and $\Omega$, that conjugates the $\RM^d$-action of their associated dynamical systems (see \cite{BBG06} section 2).

\begin{defini}
\label{Ktiling07.def-Bp} 
The {\em patch space} \(\Bb_p\) is the collared prototile space of $\Tt_p$ (definition \ref{Ktiling07.def-B0}):
\(
\Bb_p = \Bb_0^c (\Tt_p) \,.
\)
\end{defini}

\noindent The images in $\Bb_p$ of the patches $p_j$'s (tiles of $\Tt_p$) will be denoted $\pi_j$ and still be called patches.
The map \( \pG_{0,\Tt_p}^c : \Omega \rightarrow \Bb_p\), built in proposition \ref{Ktiling07.prop-projhullB0},
is denoted $\pG_p$.

\vspace{.1cm}

\noindent This construction of $\Bb_p$ from $\Tt_p$ is essentially the same as the construction of $\Bb_0$ from $\Tt$ given in definition
\ref{Ktiling07.def-B0}, the only difference being that collared prototiles of $\Tt_p$ (collared patches of $\Tt$)  are used here instead. Such
patch spaces are said to {\it force their borders}. This condition was introduced by Kellendonk in \cite{Ke95} for substitution tilings and was required
in order to be able to recover the hull as the inverse limit of such spaces. It was generalized in \cite{BBG06} for branched manifolds of
repetitive tilings with FLC and coincides here with the above definition.

\vspace{.1cm}

\noindent The map \( f_p : \Bb_p \rightarrow \Bb_0 \) defined by \( f_p = \pG_0 \circ \pG_{p}^{-1} \) is surjective and continuous.
It projects $\Bb_p$ onto $\Bb_0$ in the obvious way: a point $x$ in  $\Bb_p$ belongs to some patch $\pi_j$, hence to some tile,
and $f_p$ sends $x$ to the corresponding point in the corresponding tile $\tau_{j'}$. More precisely, if $\tilde{x}$ in $p_j$ is
the point of $\RM^d$ corresponding to $x$ in $\pi_j$, then $\tilde{x}$ belongs to some tile, which is the translate of some $t_{j'}$
and $f_p(x)$ is the corresponding point in $\tau_{j'}$ in $\Bb_0$. If $x$ belongs to the boundaries of say \(\pi_{j_1}, \cdots \pi_{j_k}\),
in $\Bb_p$, then there are corresponding points $\tilde{x}_{j_1}, \cdots \tilde{x}_{j_k}\), in \(p_{j_1}, \cdots p_{j_k}\),
which are then on the boundaries of the copies of some tiles $t_{{j'}_1}, \cdots t_{{j'}_{k'}}\). The boundaries of those tiles are identified
by the map $\rho_0$ in the definition of $\pG_0$ and $f_p(x)$ is the corresponding point on the common boundaries of the
\(\tau_{{j'}_1}, \cdots \tau_{{j'}_{k'}}\).

\vspace{.1cm}

\noindent Recall the convention stated after definition \ref{Ktiling07.def-patch}: if $\hat{p}$ is a pattern, then $p$ denotes its representative
that has its puncture at the origin. Given two patterns $\hat{p}$ and $\hat{q}$ with \( q \subset p\), the map \( f_{qp} : \Bb_p \rightarrow \Bb_q \)
defined by \( f_{qp} = f_{q}^{-1} \circ f_p = \pG_q \circ \pG_{p}^{-1}\) is continuous and surjective. Given three patterns $\hat{p}, \hat{q}$
and $\hat{r}$ with \( r \subset q \subset p\) the following composition rule holds: \( f_{rq} \circ f_{qp} = f_{rp} \). Moreover given two
arbitrary patterns $\hat{q}$ and $\hat{r}$, there exists another pattern $\hat{p}$ such that $p$ contains $q$ and $r$ (since \(\Ll_q\) and
\(\Ll_r\) are repetitive sub-Delone sets of $\Ll_\Tt$).  Hence the index set of the maps $f_p$'s is a directed set and \( ( \Bb_p , f_{qp} ) \) is a
projective system.

\vspace{.1cm}

\noindent As shown in \cite{BBG06}, the hull $\Omega$ can be recovered from the inverse limit \( \varprojlim ( \Bb_p , f_{qp} ) \) under some technical conditions
that are for convenience taken here to be directly analoguous to those given in \cite{BBG06} section 2.6.
Namely: the patch spaces $\Bb_p$ are required to {\it force their borders}, and only maps $f_{qp}$ between patch spaces 
that are {\it zoomed out} of each other are allowed (definition \ref{Ktiling07.def-zoomout} below). The first condition, as mentioned above,
is fullfilled here by the very definition of patch spaces given above, because they are built out of {\it collared patterns}
as can be checked from the more general definition of a branched manifold that forces its border given in \cite{BBG06} definition 2.43.
The second condition is stated as follows.

\begin{defini}
\label{Ktiling07.def-zoomout}
Given two patterns \(\hat{p},\hat{q} \) in \(\Pp_\Tt\), with \( q \subset p\),  \( \Bb_p \) is said to be {\rm zoomed out} of \( \Bb_q \)
if the following two conditions hold.
\begin{enumerate}[(i)]

\item For all \( i \in \{ 1, \cdots N_p\} \), the patch $p_i$ is the union of some copies of the patches $q_j$'s.

\item For all \( i \in \{ 1, \cdots N_p\} \), the patch $p_i$ contains in its interior a copy of some patch $q_j$.

\end{enumerate}
\end{defini}

\noindent The first condition is equivalent to requiring that the tiles of $\Ll_p$ are patches of $\Ll_q$. 

\vspace{.1cm}

\noindent Given a patch space $\Bb_q$ it is always possible to build another patch space $\Bb_p$ that is zoomed out of $\Bb_q$: it suffices to choose
a pattern $\hat{p}$ such that $p \supset q$ with a radius large enough. This can be done by induction for instance: building $\Bb_p$ from patches of
$\Ll_q$ (which are patches of $\Tt$) the same way $\Bb_q$ was built from patches of $\Tt$; if the radius of $q$ is large enough, then each patch of $\Tt$ that
$\Bb_q$ is composed of will contain a tile of $\Tt$ in its interior and so $\Bb_q$ will be zoomed out from $\Bb_0$.

\vspace{.1cm}

\begin{defini}
\label{Ktiling07.def-properseq}
A {\rm proper sequence} of patch spaces of $\Tt$ is a projective sequence \( \bigl\{ \Bb_{l}, f_l \bigr\}_{l\in \NM} \) where, for all $l\ge1$, $\Bb_l$ is a patch 
space associated with a pattern $\hat{p}_l$ of $\Tt$ and  \(f_l= f_{p_{l-1}p_{l}}\), such that $\Bb_l$ is zoomed out of $\Bb_{l-1}$, with the
convention that \(f_0 = f_{p_1}\) and $\Bb_0$ is the prototile space of $\Tt$.
\end{defini}

\noindent Note that the first patch space in a proper sequence can be chosen to be the prototile space of the tiling ({\it i.e.} made of uncollared prototiles),
all that matters for recovering the hull by inverse limit as shown in the next theorem, is that the next patch spaces are zoomed out of each other, and built
out of collared patches ({\it i.e.} force their borders).

\begin{theo}
\label{Ktiling07.thm-invlimBp}
The inverse limit of a proper sequence \( \bigl\{ \Bb_{l}, f_l \bigr\}_{l\in \NM} \) of patch spaces of $\Tt$ is homeomorphic to the hull of $\Tt$: 
\[
\Omega \cong \varprojlim \bigl( \Bb_{l}, f_l \bigr)  \,.
\]
\end{theo}

\noindent {\it Proof.} The homeomorphism is given by the map \( \pG : \Omega \rightarrow \varprojlim   \bigl( \Bb_{l}, f_l \bigr) \), defined by
\( \pG (\omega) = \bigl( \pG_0 (\omega), \pG_1 (\omega), \cdots \bigr) \), with inverse 
\( \pG^{-1} (x_0, x_1, \cdots ) = \cap \{ \pG^{-1}_l (x_l), l\in \NM\} \).

\noindent The map $\pG$ is surjective, because all the $\pG_l$'s are. For the proof of injectivity, consider \( \omega, \omega' \in \Omega\) with
\( \pG(\omega) = \pG(\omega') \). For simplicity the metric $\delta$ defined in remark \ref{Ktiling07.rem-topohull} is used in this proof.
Fix $\epsilon >0$. To prove that $\delta(\omega, \omega') < \epsilon$ it suffices to show
that the two tilings agree on a ball of radius $\frac{1}{\epsilon}$.
For each $l$ in $\NM$, \(\pG_l (\omega) = \pG_l (\omega') \) in some patch $\pi_{l,j}$ in $\Bb_l$. This means that the tilings agree on some 
translate of the patch $p_{l,j}$ that contains the origin. The definition of patch spaces made up from collared patterns (condition of {\it forcing the border}
in \cite{BBG06}) implies that the tilings agree on the ball $B(0_{\RM^d}, r_l)$, where $r_l$ is the parameter of uniform discretness of the Delone set
$\Tt_{p_l}$. The assumption that the patch spaces are zoomed out of each other (condition (ii) in definition \ref{Ktiling07.def-zoomout}) implies that 
$r_l > r_{l-1} + r$, where $r$ is the parameter of uniform discretness of $\Ll_\Tt$. Hence $r_l > (l+1) r$, and choosing $l$ bigger than the integer part of 
$\frac{1}{r\epsilon}$ concludes the proof of the injectivity of $\pG$.

\noindent Condition (i) in definition \ref{Ktiling07.def-zoomout} implies that for every $l\ge 1$, \( \pG^{-1}_l (x_l) \subset \pG^{-1}_{l-1} (x_{l-1})\).
The definition of $\pG^{-1}$ then makes sense by compactness of $\Omega$ because any finite intersection of some $\pG^{-1}_l (x_l)$'s is
non empty and closed. 

\noindent A neighborhood of $x=(x_0, x_1, \cdots)$ in  \(\varprojlim   \bigl( \Bb_{l}, f_l \bigr)\) is given by
\( U_n (x) = \bigl\{ y=(y_0, y_1, \cdots) \, : \,  y_i = x_i , i \le n \bigr\} \) for some integer $n$.  If $\omega$ is a preimage of $x$, \( \pG(\omega) = x\),
then the preimage of $U_n(x)$ is given by all tilings $\omega'$ such that \(\pG_n (\omega') = \pG_n (\omega)\), {\it i.e.} those tilings agree with 
$\omega$ on some patch $p_n$ around the origin; they form then a neighborhood of $\omega$ in $\Omega$. This proves that
that $\pG$ and $\pG^{-1}$ are continuous. \qed

\vspace{.1cm}

\noindent Another important construction which allows one to recover the hull by inverse limit has been given by G\"ahler in an unpublished work
as a generalization of the construction of Anderson and Putnam in \cite{AP98}. Instead of gluing together patches to form the patch spaces $\Bb_p$'s,
G\"ahler's construction consists in considering ``multicollared prototiles'' and the graph that link them according to the local configurations of their
representatives in the tiling. This construction keeps track of the combinatorics of the patches that those ``multicollared tiles''
represent and thus is enough to recover the hull topologically, which is sufficient for topological concerns (for cohomology or $K$-theory in particular).
The construction given here, and taken from \cite{BBG06} where those patch spaces are proven to be branched manifolds, takes more structure into
account and the homeomorphism between the hull and the inverse limit of such branched manifolds which is built in \cite{BBG06} is not only a
topological conjugacy but yields a conjugacy of the dynamical systems' actions.

\section{Cohomology of Tiling Spaces}
\label{Ktiling07.sec-cohomtil}

\subsection{Tiling Cohomologies}
\label{Ktiling07.ssect-cohom}

\noindent Various cohomologies for tiling spaces have been used in the literature. For the class of tilings considered here (aperiodic, repetitive with FLC)
they are all isomorphic.

\vspace{.1cm}

\noindent First, the \v{C}ech cohomology of the hull was introduced for
instance in \cite{AP98} for substitution tilings (and in full generality for tilings on Riemanian manifolds in \cite{Sa03}).
Hulls of such tilings are obtained by inverse limit of finite $CW$-complexes, and their \v{C}ech cohomology is obtained by direct limit.

\vspace{.1cm}

\noindent For repetitive tilings with FLC,  if a lamination structure is given to the hull as in \cite{BBG06}, the cohomology of the hull is defined by
direct limit of the simplicial cohomologies of branched manifolds that approximate the hull by inverse limit (a generalisation of theorem
\ref{Ktiling07.thm-invlimBp}). This cohomology is isomorphic the \v{C}ech cohomology of the hull, using the natural isomorphism between
\v{C}ech and simplicial cohomologies that holds for such branched manifolds (which are $CW$-complexes) and passing to direct limit.
It has been used to prove that the generators of the $d$-th cohomology group are in one-to-one correspondence with invariant ergodic probability
measures on the hull.

\vspace{.1cm}

\noindent Another useful cohomology, the Pattern-Equivariant (PE) cohomology, has been proposed by Kellendonk and Putnam in \cite{Ke03,KP06} for
real coefficients and then generalized to integer coefficients by Sadun \cite{Sa06}. This cohomology has been used for proving that the Ruelle-Sullivan map
(associated with an ergodic invariant probability measure on the hull) from the \v{C}ech cohomology of the hull to the exterior algebra of the dual of $\RM^d$
is a ring homomorphism. Let $\Tt$ be a repetitive tiling with FLC. Assume its tiles
are compatible $CW$-complexes, so that it gives a $CW$-complex decomposition of $\RM^d$. The group of PE $n$-cochains $C_{PE}^n$ is a subgroup
of the group of integer singular $n$-cochains of $\RM^d$ that satisfy the following property:
an $n$-cochain $\varphi$ is said to be PE if there exists a patch $p$ of $\Tt$ such that
\(\varphi(\sigma_1) = \varphi(\sigma_2)\), for $2$ $n$-simplices \(\sigma_1, \sigma_2\) with image cells $e_1, e_2$,
whenever there exists \( x \in e_1, \, y \in e_2\) such that \(\pG_p ( \tra^{-x}\Tt ) = \pG_p ( \tra^{-y} \Tt )\).
The simplicial coboundary of a PE cochain is easily seen to be PE (possibly with respect to a patch of larger radius), and this defines the complex
for integer PE cohomology. A PE cochain is a cochain that agree on points which have the same local environments in $\Tt$ or equivalently 
is the pull back of a cochain on some patch space $\Bb_p$. Hence PE cohomology can be seen as the direct limit of the singular cohomologies of
a proper sequence of patch spaces. Using the natural isomorphisms between cellular and \v{C}ech cohomologies that holds for those $CW$-complexes and
taking direct limits, the integer PE cohomology turns out to be isomorphic to the \v{C}ech cohomology of the hull.

\vspace{.1cm}

\noindent The PV cohomology described in the next section \ref{Ktiling07.ssect-ringtrans} is also isomorphic to the \v{C}ech cohomology of the hull
(see theorem \ref{Ktiling07.thm-pimcoho}).

\subsection{The PV cohomology}
\label{Ktiling07.ssect-ringtrans}

\noindent The definition of a {\it $\Delta$-complex} structure is first recalled, following the presentation of Hatcher in his book on Algebraic Topology \cite{Hatcher02} section 2.1.

\vspace{.1cm}

\noindent Given $n+1$ points $v_0, \cdots v_n$, in $\RM^m, m>n$, that are not collinear, let $[v_0, \cdots, v_n]$ denote the $n$-simplex with vertices
$v_0, \cdots v_n$. Let $\Delta^n$ denote the standard $n$-simplex

\[
\Delta^n = \bigl\{  (x_0, x_1, \cdots x_n) \in \RM^{n+1} \, : \, \sum_{i=0}^{n} x_i = 1 \ {\rm and} \ x_i \ge 0 \ {\rm for \ all} \ i \bigr\}\,,
\]

\noindent whose vertices are the unit vectors along the coordinate axes. An ordering of those vertices is specified and this allows to 
define a canonical linear homeomorphism between $\Delta^n$ and any other
$n$-simplex $[v_0, \cdots, v_n]$, that preserves the order of the vertices, namely, \(  (x_0, x_1, \cdots x_n) \mapsto \sum x_i v_i \).

\vspace{.1cm}

\noindent If one of the $n+1$ vertices of an $n$-simplex $[v_0, \cdots, v_n]$ is deleted, then the remaining $n$ vertices span an
$(n-1)$-simplex, called a face of $[v_0, \cdots, v_n]$. By convention the vertices of any subsimplex spanned by 
a subset of the vertices are ordered according to their order in the larger simplex.

\vspace{.1cm}

\noindent The union of all the faces of $\Delta^n$ is the boundary of $\Delta$, written $\partial \Delta^n$. The open simplex
\( \stackrel{\! \! \! \! \circ}{\Delta^n} \) is \( \Delta^n \setminus \partial \Delta^n \), the interior of $\Delta^n$.

\vspace{.1cm}

\noindent A $\Delta$-complex structure on a space $X$ is a collection of maps \( \sigma_\alpha : \Delta^n \rightarrow X\), with
$n$ depending on the index $\alpha$, such that

\begin{enumerate}[(i)]
\item The restriction \( \sigma_\alpha \arrowvert_{\stackrel{\! \! \! \! \circ}{\Delta^n}} \) is injective, and each point of $X$ is in the image of exactly 
one such restriction.
\item  Each restriction of $\sigma_\alpha$ to a face of $\Delta^n$ is one of the maps \( \sigma_\beta : \Delta^{n-1} \to X\). The face of
$\Delta^n$ is identified with $\Delta^{n-1}$ by the canonical linear homeomorphism between them that preserves the ordering of
the vertices.
\item A set $A \subset X$ is open iff \( \sigma_{\alpha}^{-1} (A)\) is open in $\Delta^n$ for each $\sigma_\alpha$.
\end{enumerate}

\noindent A $\Delta$-complex $X$ can be built as a quotient space of a collection of disjoint simplices by identifying various subsimplices spanned 
by subsets of the vertices, where the identifications are performed using the canonical linear homeomorphism that preserves orderings of
the simplices.  It can be shown that  $X$ is a Hausdorff space. Then condition (iii) implies that each restriction
\( \sigma_\alpha \arrowvert_{\stackrel{\! \! \! \! \circ}{\Delta^n}}  \) is a homeomorphism onto its image, which is thus an open simplex in 
$X$. These open simplices \( \sigma_\alpha ( \stackrel{\! \! \! \! \circ}{\Delta^n} ) \) are cells $e^{n}_{\alpha}$ 
of a $CW$-complex structure on $X$ with the $\sigma_\alpha$'s as characteristic maps.

\vspace{.1cm}

\noindent  It is now assumed that the tiles of $\Tt$ are compatible finite $\Delta$-complexes: the intersection of two tiles is a sub-$\Delta$-complex
of both. In addition each cell \(e^{n}_{\alpha}=\sigma_\alpha ( \stackrel{\! \! \! \! \circ}{\Delta^n} ) \) of each tile is punctured, by say the image
under $\sigma_\alpha$ of the barycenter of $\Delta^n$. Hence $\Tt$ can be seen as a $\Delta$-complex decomposition of $\RM^d$, and this
$\Delta$-complex structure gives a ``refinement'' of the tiling (each tile beeing decomposed into the union of the closures of the cells it contains).
If $\Tt$ is the Voronoi tiling of a Delone set its tiles can be split into $d$-simplices since they are are polytopes:
this gives a $\Delta$-complex structures for which the maps $\sigma_\alpha$ are simply affine maps. The maps \(\sigma_\alpha : \Delta^n \rightarrow \Bb_0\)
of the $\Delta$-complex structure of $\Bb_0$, will be called the {\em characteristic maps} of the $n$-simplices on $\Bb_0$ represented by its image.

\vspace{.1cm}

\noindent By construction (definition \ref{Ktiling07.def-B0}) as a quotient $CW$-complex, the prototile space $\Bb_0$ is a finite $\Delta$-complex.

\begin{defini}
\label{Ktiling07.def-simplicialtrans}
The {\em $\Delta$-transversal}, denoted $\Xi_\Delta$, is the subset of $\Omega$ consisting of tilings that have the puncture
of one of their cells (of one of their tiles) at the origin $0_{\RM^d}$.
\end{defini}

\noindent Since $\Tt$ has finitely many prototiles and they have finite $CW$-complex structures, the set of the punctures of the cells of the tiles
of $\Tt$ is a Delone set. This Delone set will be called the {\it $\Delta$-Delone set} of $\Tt$, and denoted $\Ll_\Delta$.
The $\Delta$-transversal is thus the canonical transversal of $\Ll_\Delta$.

\vspace{.1cm}

\noindent The $\Delta$-transversal is not immediately related to the canonical transversal. Indeed the lift of
the puncture of a cell does not belong to the transversal in general, unless this puncture coincides with the puncture of the tile
of $\Bb_0$ that contains that cell.

\vspace{.1cm}

\noindent The $\Delta$-transversal is the lift of the punctures of the cells in $\Bb_0$. It is partitioned by the lift of the punctures
of the $n$-cells, denoted $\Xi_\Delta^n$ which is the subset of $\Omega$ consisting of tilings that have the puncture of an $n$-cell at the origin.
As for the transversal, the $\Delta$-transversal is a Cantor set, and the $\Xi_\Delta^n$'s give a partition by clopen sets.
The ring of continuous integer-valued functions on the $\Delta$-transversal, \(C (\Xi_\Delta, \ZM)\), is thus the direct sum of the
\(C (\Xi_\Delta^n, \ZM)\)'s for $n=0, \cdots d$.

\vspace{.1cm}

\noindent Given the characteristic map $\sigma$ of an $n$-simplex $e$ of $\Bb_0$, let $\Xi_\Delta (\sigma)$ denote the lift of the puncture of $e$, and
$\chi_\sigma$ its characteristic function on $\Xi_\Delta$.  The subset $\Xi_\Delta (\sigma)$ is called the {\it acceptance zone} of $\sigma$.
Since continuous integer-valued functions on a totally disconnected space are generated by characteristic
functions of clopen sets, $\chi_\sigma$ belongs to \(C (\Xi_\Delta, \ZM)\).

\vspace{.1cm}

\noindent Consider the characteristic map $\sigma: \Delta^n \rightarrow \Bb_0$ of an $n$-simplex $e$ of $\Bb_0$, and $\tau$ a face of $\sigma$ ({\it i.e.} the restriction of
$\sigma$ to a face of $\Delta^n$) with associated simplex $f$ (a face of $e$).
The simplices $e$ and $f$ in $\Bb_0$ are contained in some tile $\tau_j$. Viewing $e$ and $f$ as subsets of the tile $t_j$ in $\RM^d$,  it is possible to
define the vector $x_{\sigma \tau}$ that joins the puncture of $f$ to the puncture of $e$. Notice that since $\Bb_0$ is a {\it flat} branched manifold
(remark \ref{Ktiling07.rem-lamination}, and \cite{BBG06}), the vector $x_{\sigma \tau}$ is also well defined in $\Bb_0$ as a vector in a region containing the simplex $e$.

\begin{defini}
\label{Ktiling07.def-ringtrans}
(i) Let $\sigma$ and $\tau$ be characteristic maps of simplices in $\Bb_0$. The operator \(\theta_{\sigma \tau}\), on \(C (\Xi_\Delta, \ZM)\), is defined by
\[
\theta_{\sigma \tau} = 
\left\{ 
\begin{array}{clc}
\chi_\sigma  \tra^{  x_{\sigma \tau} } \chi_{\tau}    & {\rm if} & \tau \subset \partial \sigma \,,\\
0 & {\rm otherwise} \,. &
\end{array}
\right.
\]
where $\tau \subset \partial \sigma$ means that $\tau$ is a face of $\sigma$. Here the translation acts by $\tra^{x_{\sigma \tau}}f(\xi) = f(\tra^{-x_{\sigma \tau}}\xi)$
whenever $f\in C (\Xi_{\Delta}(\tau), \ZM)$.

\noindent (ii) The {\rm function ring of the transversal} $\Aa_{\Xi_\Delta}$ is the ring (finitely) generated by
the operators \(\theta_{\sigma \tau}\) and their adjoints \(\theta_{\sigma \tau}^\ast = \chi_\tau \tra^{-x_{\sigma\tau}} \chi_\sigma\) if \( \tau \subset \partial \sigma\) and $0$ otherwise,
over all characteristic maps $\sigma$ and $\tau$ of simplices in $\Bb_0$.
\end{defini}

\noindent The operators $\theta_{\sigma \tau}$'s satisfy the following properties:

\begin{subequations}
\label{Ktiling07.eq-ringtrans}
\begin{align}
\theta_{\sigma \tau} \theta_{\sigma\tau}^\ast =  \chi_\sigma \quad {\rm if} \ \tau \subset \partial \sigma \,,
\label{Ktiling07.eq-ringtrans2}\\
\sum_{\sigma \, : \, \partial \sigma \supset \tau}   \theta_{\sigma\tau}^\ast \theta_{\sigma \tau} =    \chi_\tau \,. \quad
\label{Ktiling07.eq-ringtrans3}
\end{align}
\end{subequations}

\noindent The $\theta_{\sigma \tau}$'s are thus {\it partial isometries}. The ring $\Aa_{\Xi_\Delta}$ is unital and its unit
is the characteristic function of the $\Delta$-transversal \( 1_{\Aa_{\Xi_\Delta}} = \chi_{\Xi_\Delta} \).
The ring of continuous integer-valued functions on the $\Delta$-transversal  \(C (\Xi_\Delta, \ZM)\) is a left $\Aa_{\Xi_\Delta}$--module.

\begin{rem}
\label{Ktiling07.rem-ringtrans}
{\em It is important to notice that the operators \(\theta_{\!\sigma \partial_i \!\sigma}\)'s are in one-to-one correspondence with a set of
arrows that generate $\Gamma_\Delta^1$, the arrows of the groupoid $\Gamma_\Delta$ of the $\Delta$-transversal (groupoid of the transversal of the Delone
set $\Ll_\Delta$, see definition \ref{Ktiling07.def-groupoid}). The set of arrows $\Gamma_\Delta^1$ is indeed in one-to-one correspondence with the set of
vectors joining points of $\Ll_\Delta$, and it is generated by the set of vectors joining a point of $\Ll_\Delta$ to a ``nearest neighbour''
(a puncture of a face of the simplex of $\Tt$ whose puncture is that point, or of a simplex containing it on its boundary).
Each vector can be decomposed into a sum of such generating vectors, and
each arrow can be decomposed into a composition of such generating arrows (corresponding to those generating vectors).
This generating set is finite, since $\Tt$ has FLC (in particular there are finitely
many prototiles), and each of its vectors is an $x_{\!\sigma \partial_i \!\sigma}$ and thus corresponds to a unique operator
\(\theta_{\!\sigma \partial_i \!\sigma}\).

\noindent A ``representation by partial isometries'' of the set of generators of $\Gamma_\Delta^1$ is given by
\[ 
\theta(\gamma) = \chi_{\pG_0 \circ \pG_0^{-1} (s(\gamma))} \, \tra(\gamma) \, \chi_{\pG_0 \circ \pG_0^{-1} (r(\gamma))} \,,
\]
\noindent where, if \( \gamma=(\xi, x)\), then $\tra(\gamma)$ is the translation operator $\tra^x$.}
\end{rem}

\noindent Let $\Ss_0^n$ be the set of the characteristic maps \(\sigma : \Delta^n \rightarrow \Bb_0\) of the $n$-simplices of the $\Delta$-complex decomposition of $\Bb_0$,
and $\Ss_0$ the (disjoint) union of the $\Ss_0^n$'s. The group of simplicial $n$-chains on $\Bb_0$, \(C_{0,n}\), is the free abelian group with basis $\Ss_0^n$.

\begin{defini}
\label{Ktiling07.def-pimcoho}
The {\rm PV cohomology} of the hull of $\Tt$ is the homology of the complex \( \bigl\{ C^n_{PV}, d_{PV}^n \bigr\} \), where:

\begin{enumerate}[(i)]

\item the PV cochain groups are the groups of continuous integer valued functions on $\Xi_\Delta^n$:
\( C^n_{PV} = C (\Xi_\Delta^n, \ZM)\) for $n=0, \cdots d$,

\item the PV differential, $d_{PV}$,  is the element of $\Aa_{\Xi_\Delta}$ given by the sum over  $n=1, \cdots d$, 
of the operators
\begin{equation}
\label{Ktiling07.eq-pimdif}
d_{PV}^n : 
\left\{ 
\begin{array}{l}
C^{n-1}_{PV} \longrightarrow C^n_{PV} \\
d_{PV}^n = {\displaystyle \sum_{\sigma \in \Ss_{0}^n} \sum_{i=0}^{n}} (-1)^i \; \theta_{\!\sigma \partial_i\!\sigma}
\end{array}
\right. \,.
\end{equation}

\end{enumerate}
\end{defini}

\noindent The ``simplicial form'' of $d_{PV}^n$ makes it clear that \(d_{PV}^{n+1} \circ d_{PV}^{n} = 0\) for $n=1, \cdots d-1$. We shall call alternatively the PV cohomology of the hull of $\Tt$ simply
the {\it PV cohomology of $\Tt$}.

\begin{rem}
\label{Ktiling07.rem-fibrahull}
{\em Thanks to the lamination structure on $\Omega$ described in remark \ref{Ktiling07.rem-lamination} the map $\pG_0$

\[
\begin{array}{ccl}
\Xi_\Delta & \hookrightarrow & \Omega \\
 & & \downarrow  \pG_0 \\
& & \Bb_0
\end{array}
\]

\noindent looks very much like a fibration of $\Omega$ with base space $\Bb_0$ and fiber $\Xi_\Delta$. However because of the branched structure there are paths in $\Bb_0$
that cannot be lifted to any leaf of $\Omega$. Hence this is not a fibration. Nevertheless the result in the present paper (theorem \ref{Ktiling07.thm-SSpim}) gives a spectral sequence analogous
to the Serre spectral sequence for a fibration. 
}
\end{rem}

\begin{rem}
\label{Ktiling07.rem-localcoef}
{\em This PV cohomology of $\Tt$ is formally analogous to a cohomology with local coefficients (see \cite{Hatcher02} chapter 3.H.) but in a more general setting.
Let $G_x$ be the group \(C( \pG_0^{-1}(x), \ZM)\) for $x$ in $\Bb_0$. As $\pG_0^{-1}(x)$ is a Cantor set, $G_x$ is actually its
$K^0$-group: \( G_x = K^0 \bigl( \pG_0^{-1}(x) \bigr)\) while its $K^1$ group is trivial. The family of groups \( (G_x)_{x\in\Bb_0}\) is analoguous to a local
coefficient system, with the operators $\theta_{\sigma\tau}$'s in the role of group isomorphisms between the fibers, but they
are not isomorphisms here (as they come from a groupoid). If $\sigma$ in $C_0^n$ is an $n$-chain, and $\varphi$ in $C_{PV}^{n-1}$ is a
PV $(n\!-\!1)$-cochain, then the differential
\[
d_{PV}^n \varphi (\sigma) =  \sum_{i=0}^{n} (-1)^i \theta_{\!\sigma \partial_i\!\sigma} \varphi (\partial_i \sigma) \,,
\]

\noindent takes the restrictions of $\varphi$ to the faces $\partial_i \sigma$, which are elements of the groups $G_{x_{\partial_i \sigma}}$
(where $x_{\partial_i \sigma}$ is the puncture of the image simplex of $\partial_i \sigma$), and pull them to $\sigma$ with the operators
\(\theta_{\!\sigma \partial_i\!\sigma}\)'s to get an element of the group $G_{x_{\sigma}}$
(where $x_{\sigma}$ is the puncture of the image simplex of $\sigma$).

\noindent A rigorous formulation of the points in this remark as well as in the previous remark \ref{Ktiling07.rem-ringtrans}
will be investigated in further research.}
\end{rem}

\subsection{An example: PV cohomology in $1$-dimension}
\label{Ktiling07.ssect-PimCoho1d}

\noindent This section presents the explicit calculation of PV cohomology of some $1$-dimensional tilings obtained by ``cut and projection''. This is to illustrate techniques that can be used to calculate PV cohomology, rather than showing its distinct features and differences compared to other cohomologies. Further examples in higher dimensions, and methods of calculation for PV cohomology will be investigated in a future work.


\begin{figure}[!h]
\psfrag{A}{$L_\alpha^\perp$}
\psfrag{B}{$L_\alpha$}
\psfrag{C}{$W_\alpha$}
\psfrag{D}{$C$}
  \begin{center}
  \includegraphics[width=7cm]{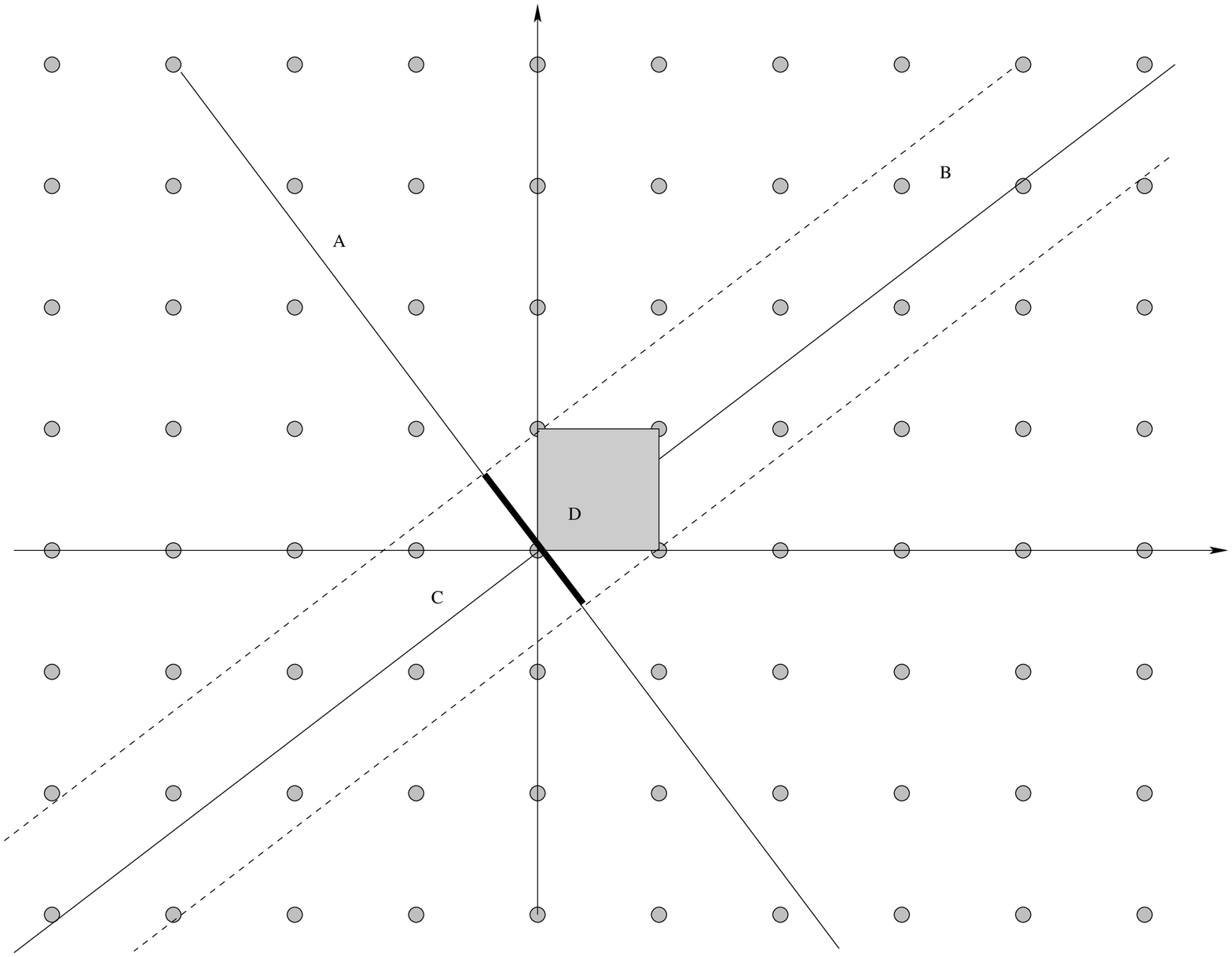}
  \end{center}
\end{figure}

\noindent Let $\RM^2$ be the Euclidean plane with basis vectors $e_1$ and $e_2$, and let $\ZM^2$ be the lattice of points of integer coordinates.
Let \( C = (0,1] \times [0,1)\) be the unit square ``open on the left and top''.
Fix an irrational number $\alpha$ in \(\RM_+ \setminus \QM\).
Let $L_\alpha$ be line of slope \(\alpha\) (through the origin) and $p_\alpha$ the projection onto $L_\alpha$.
Let $L_\alpha^\perp$ be the orthogonal complement of $L_\alpha$ and \(p_\alpha^\perp = 1 - p_\alpha\).
Let $W_\alpha=p^{\perp}_\alpha ( C )$, and \( \Sigma = \bigl\{ a \in \ZM^2 \, : \, p_\alpha^\perp(a) \in W_\alpha \bigr\}\) be the set of points of $\ZM^2$ that project to $W_\alpha$.
A lattice point $a = (a_1, a_2)$ in $\ZM^2$ belongs to $\Sigma$ if and only if \( -\alpha \le a_2 - a_1 \alpha  < 1 \). 
The projection of $\Sigma$ onto $L_\alpha$ defines an aperiodic tiling $\Tt_\alpha$ whose tiles' edges are the projections of points of $\Sigma$. 
The tiling $\Tt_\alpha$ has two types of tiles that correspond to the projections of the vertical and horizontal of $C$.
In particular $\Tt_\alpha$ is repetitive with FLC \cite{BHZ00,La99A, La99B, LP03}. 

\vspace{.1cm}

\noindent Identifying $L_\alpha^\perp$ (oriented towards the ``northwest'') with the real line, $W_\alpha$ can be seen as the interval $[-\frac{\alpha}{\sqrt{1+\alpha^2}}, \frac{1}{\sqrt{1+\alpha^2}})$.
Let $\Aa_\alpha$ be the \Cs generated by the characteristic functions $\chi_{[x_1,x_2)}$ for $x_1<x_2$ in $p_\alpha^\perp(\Sigma)$.
The transversal $\Xi_\alpha$ is defined to be the spectrum of $\Aa_\alpha$.
It has been shown in \cite{BHZ00} that it is a Cantor set (compact, perfect and totally disconnected), and can be seen as a completion of the interval $W_\alpha$ for a finer topology than the usual one, where the intervals $[x_1,x_2)$, for $x_1<x_2$ in $p_\alpha^\perp(\Sigma)$, are open and closed.

\vspace{.1cm}

\noindent We further identify $W_\alpha$ (rescalling it and identifying its endpoints) with the unit circle with the topology induced by $\Xi_\alpha$.
This defines the Cantor circle $\SM^1_\alpha$.
It admits a countable basis of open and closed sets \(i_{lm} = [ l \frac{\alpha}{1+\alpha}, m \frac{\alpha}{1+\alpha}) \mod 1\) for integers $l,m$.
We can alternatively define $\SM^1_\alpha$ as the spectrum of the \Cs generated by $\{ \theta^n_\alpha \chi_\alpha, \, n\in \ZM\}$ where $\chi_\alpha$ is the characteristic function of the arc $[0,\frac{\alpha}{1+\alpha})$, and $\theta_\alpha$ the rotation of angle $\frac{\alpha}{1+\alpha}$.

\vspace{.1cm}

\noindent We give $L_\alpha$ an orientation towards the first quadrant (``northeast''), and puncture the tiles by their left vertices for convenience (not their barycenters).
The prototile space $\Bb_0 = \SM^1 \!\vee \SM^1$ is the wedge sum of two circles corresponding to the two prototiles with their vertices identified.
The difference with section \ref{Ktiling07.ssect-ringtrans} is that the $\Xi_\Delta^n$'s do not partition $\Xi_\Delta$ since there is only one vertex here, see definition \ref{Ktiling07.def-simplicialtrans}.
Consequently the canonical transversal and $\Delta$-transversal are identical here.
The transversal $\Xi_\alpha$ is seen as the Cantor circle $\SM^1_\alpha$.
The arc \( [ \frac{-\alpha}{1+\alpha}, 0)\) is the acceptance zone $\Xi_a$ of the prototile $a$ whose representatives
are the projections of horizontal intervals of $\Sigma$, and the arc \( [0, \frac{1}{1+\alpha})\) is the acceptance zone $\Xi_b$ of the prototile whose representatives are the projections of vertical intervals of $\Sigma$.
The vectors $x_{\sigma \tau}$'s in the definitions of the operators \(\theta_{\sigma\tau}\)'s are here induced by the projections \(p_\alpha^\perp (e_1)\) and \(p_\alpha^\perp (e_2)\). 
Once rescaled to the unit circle the operator $\tra^{x_a}$ becomes the rotation by \(\frac{-1}{1+\alpha}\), and $\tra^{x_b}$ the rotation by $\frac{\alpha}{1+\alpha}$, and are thus equal.
The PV complex reads here simply
\begin{equation}
\label{Ktiling07.eq-complexpim1d}
0 \longrightarrow C \bigl( \SM_\alpha^1, \ZM \bigr) \xrightarrow{d_{PV} = {\rm id} - \theta_\alpha} C \bigl( \SM_\alpha^1, \ZM \bigr) \longrightarrow 0 \,,
\end{equation}
where $\theta_\alpha$ is the rotation by $\frac{\alpha}{1+\alpha}$ on $\SM_\alpha^1$, and is unitary.

\begin{proposi}
\label{Ktiling07.prop-pimcoho1d}
The PV cohomology of (the hull $\Omega_\alpha$ of) $\Tt_\alpha$ is given by
\[
\left\{
\begin{array}{lll}
H_{PV}^0 \bigl( \Bb_0; C( \SM_\alpha^1, \ZM) \bigr) & \cong & \ZM \\
H_{PV}^1 \bigl( \Bb_0; C( \SM_\alpha^1, \ZM) \bigr) & \cong & \ZM \oplus \ZM
\end{array}
\right. \,.
\]
\end{proposi}

\noindent {\it Proof.} A function \( f \in C \bigl( \SM_\alpha^1, \ZM \bigr) \) reads as a finite sum \( f = \sum n \chi_{I_n} \) where $n$ is an integer and
$\chi_{I_n}$ is the characteristic function of the clopen set $I_n = f^{-1}(n)$. The $0$-th cohomology group is the set of invariant functions under $\theta_\alpha$.
Each $I_n$ is a finite disjoint union of base clopen sets $i_{lm}$'s.
Given two base clopens $i_{lm} \subset I_n$ and $i_{l'm'} \subset I_{n'}$,
\( \theta_\alpha^{m-m'} i_{l'm'} = [(l'+m-m')\frac{\alpha}{1+\alpha}, m\frac{\alpha}{1+\alpha}) \mod 1\), so that
\( i_{lm} \cap \theta_\alpha^{m-m'} i_{l'm'} \ne \emptyset\). Hence if $f$ is invariant under $\theta_\alpha$, then $n$ must equal $n'$ and therefore
$f$ must be constant.

\noindent The calculation of the first cohomology group, the group of coinvariants under $\theta_\alpha$, relies upon the fact that any function
\( f \in C \bigl( \SM_\alpha^1, \ZM \bigr) \) can be written as \( n_f \chi_{i_{01}} + m_f \chi_{\SM_\alpha^1}\) modulo
\( (1-\theta_\alpha) C \bigl( \SM_\alpha^1, \ZM \bigr) \), for some integers $n_f, m_f$, that are uniquely determined by the class of $f$.
This technical point is tedious but elementary (it uses an encoding of the
real numbers from the partial fraction decomposition of $\frac{\alpha}{1+\alpha}$). \qed

\begin{rem}
\label{Ktiling07.rem-pimcoho1d}
{\em
\begin{enumerate}[(i)]

\item Using integration, $\int_{\SM_\alpha^1}$, one gets a group isomorphism
\[
H_{PV}^1 \bigl( \Bb_0; C( \SM_\alpha^1, \ZM) \bigr) \cong  \ZM \oplus \frac{\alpha}{1+\alpha} \ZM \,,
\]
but PV cohomology by itself cannot distinguish between different values of $\alpha$.

\item The above calculations yield also the $K$-groups of the tiling space:
\[
\left\{
\begin{array}{c}
K^0 (\Omega) \cong \check{H}^0 (\Omega_\alpha ; \ZM) \cong H_{PV}^0 \bigl( \Bb_0; C( \SM_\alpha^1, \ZM) \bigr) \\
K^1 (\Omega) \cong \check{H}^1 (\Omega_\alpha ; \ZM) \cong H_{PV}^1 \bigl( \Bb_0; C( \SM_\alpha^1, \ZM) \bigr) 
\end{array}
\right.
\]

\noindent The isomorphism between \v{C}ech and PV cohomologies for hull of tilings is proven in theorem \ref{Ktiling07.thm-pimcoho}.
The isomorphisms between $K$-theory and cohomology of the hull comes from the natural isomorphisms (Chern character) between $K$-theory and cohomology of $CW$-complexes of dimension less than $3$ (a proof can be found in \cite{AP98}, proposition 6.2),
and the fact that $\Omega_\alpha$ is the inverse limit of such space (theorem \ref{Ktiling07.thm-invlimBp}).

\item With the previous remark, the PV complex \eqref{Ktiling07.eq-complexpim1d} can be viewed as the Pimsner-Voiculescu exact
sequence \cite{PV80} for the $K$-theory of the \Cs \(C(\SM_\alpha^1) \rtimes_{\theta_\alpha} \ZM\):
\[
\begin{array}{ccccc}
K_0 \bigl( C(\SM_\alpha^1) \bigr) & \xrightarrow{ \ {\rm id} - \theta_\alpha \  } & K_0 \bigl( C(\SM_\alpha^1) \bigr) &
\xrightarrow{ \quad \quad } & K_0 \bigl( C(\SM_\alpha^1) \rtimes_{\theta_\alpha} \ZM \bigr) \\
\uparrow & & & & \downarrow \\
K_1 \bigl( C(\SM_\alpha^1) \rtimes_{\theta_\alpha} \ZM \bigr) & \xleftarrow{ \quad \quad } & K_1 \bigl( C(\SM_\alpha^1) \bigr) &
\xleftarrow{ \ {\rm id} - \theta_\alpha \ } & K_1 \bigl( \SM_\alpha^1 \bigr)
\end{array} \,.
\]
\noindent Indeed \(C(\SM_\alpha^1) \rtimes_{\theta_\alpha} \ZM\) is the \Cs of the groupoid of the transversal, which is Morita equivalent to the \Cs
of the hull of $\Tt_\alpha$ \cite{Be86}. Also \( K_0 \bigl( C(\SM_\alpha^1) \bigr) \cong C(\SM_\alpha^1, \ZM)\) and \( K_1 \bigl( C(\SM_\alpha^1) \bigr) \cong 0\)
\cite{Be91}.
\end{enumerate}
}
\end{rem}

\section{Proof of theorem \ref{Ktiling07.thm-SSpim}}
\label{Ktiling07.sec-proof}

\noindent  The proof of theorem \ref{Ktiling07.thm-SSpim} follows from the following two theorems that will be proved separately for convenience and clarity.

\vspace{.1cm}

\noindent It is important to notice that the first page of the spectral sequence of theorem 
\ref{Ktiling07.thm-SSpim} is the PV complex (definition \ref{Ktiling07.def-pimcoho}) and therefore its second page is given by PV cohomology. However the direct identification of the PV differential on the first page is highly technical and will not be presented here. The proof will use instead an approximation of the hull (theorem \ref{Ktiling07.thm-invlimBp}) as an inverse limit of patch spaces (definition \ref{Ktiling07.def-Bp}), and a direct limit spectral sequence for the $K$-theory of the hull (theorem \ref{Ktiling07.thm-SS}), with page-$2$ isomorphic to its \v{C}ech cohomology. Some abstract results on spectral sequences are also needed, in particular some properties of the Atiyah-Hirzebruch spectral sequence \cite{AH61}. For the convenience of the reader, those known results are grouped together and proven in a separate section, section \ref{Ktiling07.sect-SS}.

\vspace{.1cm}

\noindent Let $\Tt$ be an aperiodic and repetitive tiling of $\RM^d$ with FLC (definition \ref{Ktiling07.def-repFLC}), and assume its tiles are finite compatible $\Delta$-complexes.
Let $\Omega$ be its hull (definition \ref{Ktiling07.def-hull}) and $\Xi_\Delta$ its $\Delta$-transversal
(definition \ref{Ktiling07.def-simplicialtrans}).

\begin{theo}
\label{Ktiling07.thm-SS}
There is a spectral sequence that converges to the $K$-theory of the \Cs of the hull
\[
E_2^{rs} \Rightarrow  K_{r+s+d} \bigl( C(\Omega) \rtimes \RM^d \bigr) \,,
\]
\noindent and whose page-$2$ is isomorphic to the integer \v{C}ech cohomology of the hull
\[
E_2^{rs} \cong
\left\{
\begin{array}{ll}
\check{H}^r ( \Omega ; \ZM ) & s \; {\it even,} \\
 0 & s \; {\it odd.}
\end{array}
\right.
\]
\end{theo}

\begin{theo}
\label{Ktiling07.thm-pimcoho}
The integer \v{C}ech cohomology of the hull is isomorphic to the PV cohomology of the hull
\[
\check{H}^\ast(\Omega ; \ZM) \cong H^\ast_{PV} \bigl( \Bb_0 ; C(\Xi_\Delta, \ZM) \bigr) \,.
\]
\end{theo}

\subsection{Proof of theorem \ref{Ktiling07.thm-pimcoho}}
\label{Ktiling07.ssect-proofpimcoho}

\noindent The proof of theorem \ref{Ktiling07.thm-pimcoho} follows from propositions  \ref{Ktiling07.prop-pimBp} and \ref{Ktiling07.prop-pimB0}
below.  First a PV cohomology of $\Bb_p$, writen \(H^\ast_{PV} ( \Bb_0 ; C(\Sigma_p,\ZM) )\), is defined. It is
proven to be isomorphic to the simplicial cohomology of $\Bb_p$ in proposition \ref{Ktiling07.prop-pimBp}, and then
proposition \ref{Ktiling07.prop-pimB0} establishes that the PV cohomology of the hull is isomorphic
to the direct limit of the PV cohomologies of a proper sequence of patch spaces.

\vspace{.1cm}

\noindent Let $\Ss_p^n$ be the set of the characteristic maps \(\sigma_p : \Delta^n \rightarrow \Bb_p\) of the $n$-simplices of the $\Delta$-complex decomposition of $\Bb_p$,
and $\Ss_p$ the (disjoint) union of the $\Ss_p^n$'s.
The group of simplicial $n$-chains on $\Bb_p$, \(C_{p,n}\), is the free abelian group with basis $\Ss_p^n$.

\begin{rem}
\label{Ktiling07.rem-orientBp}
{\em The map \(f_p : \Bb_p \rightarrow \Bb_0\) preserves the orientations of the simplices (the ordering of their vertices).
}
\end{rem}

\noindent Given a simplex $\sigma_p$ on $\Bb_p$, let $\Xi_{p, \Delta} (\sigma_p)$ denote the lift of the puncture of its image in $\Bb_p$ and
$\chi_{\sigma_p}$ its characteristic function.  $\Xi_{p, \Delta} (\sigma_p)$ is called the {\it acceptance zone} of $\sigma_p$; it is a clopen set
in $\Xi_\Delta$.

\begin{lemma}
\label{Ktiling07.lem-partitionXip}
Given a simplex $\sigma$ on $\Bb_0$, its acceptance zone is partitioned by the acceptance zones of its preimages $\sigma_p$'s
on $\Bb_p$
\[
\Xi_\Delta (\sigma) = \coprod_{ \sigma_p \in f_{p \, \#}^{-1} (\sigma) } \Xi_{p, \Delta} (\sigma_p) \,.
\]
where \(f_{p\#} : C_{p,n} \rightarrow C_{0,n}\) denotes the map induced by $f_p$ on the simplicial chain groups.
\end{lemma}

\noindent {\it Proof.} The union of the \(\Xi_{p, \Delta} (\sigma_p)\)'s is equal to \(\Xi_\Delta (\sigma)\), because each lift of the image simplex of a $\sigma_p$
corresponds to the lift of the image simplex of $\sigma$ that has a given local configuration (it is contained in some patch of $\Tt$ or some tile of $\Tt_p$, see
section \ref{Ktiling07.ssect-invlim}.).
If $\sigma_p$ and $\sigma'_p$ are distinct in $f_{p \,\#}^{-1} (\sigma)$, then
\(\pG_p \bigl(  \Xi_{p, \Delta} (\sigma_p) \bigr)  \cap \pG_p \bigl(  \Xi_{p, \Delta} (\sigma'_p) \bigr) \) is empty, and thus
so is \( \Xi_{p, \Delta} (\sigma_p) \cap \Xi_{p, \Delta} (\sigma'_p) \). \qed

\noindent The maps in \(f_{p \,\#}^{-1} (\sigma)\) are thus in one-to-one correspondence with the set of atoms of the partition
of \(\Xi_\Delta (\sigma)\) by the \(\Xi_{p, \Delta} (\sigma_p)\)'s, and the union over $\sigma$ in $\Ss_0^n$ of the \(f_{p \,\#}^{-1} (\sigma)\)'s
is just $\Ss_p^n$.

\vspace{.1cm}

\noindent The simplicial $n$-cochain group $C^n_p$ is $Hom(C_{p,n} , \ZM)$,
the dual of the $n$-chain group $C_{p,n}$. It is represented faithfully on the group \(C(\Xi_\Delta,\ZM)\) of
continuous integer valued functions on the $\Delta$-tranversal by

\begin{equation}
\label{Ktiling07.eq-reppimcochain}
\rho_{p,n} : 
\left\{
\begin{array}{ccc}
C_p^n & \longrightarrow & C(\Xi_\Delta^n,\ZM) \\
\psi & \longmapsto & {\displaystyle \sum_{\sigma_p \in \Ss_p^n}} \psi (\sigma_p) \chi_{\sigma_p}
\end{array}
\right. \,.
\end{equation}

\noindent The image of $\rho_{p,n}$ will be written \(C(\Sigma_p^n, \ZM)\), to remind the reader that this consists of functions on the ``discrete transversal''
$\Sigma_p^n$ which corrresponds to the atoms of the partitions of the \(\Xi_\Delta (\sigma)\)'s by the \(\Xi_{p,\Delta} (\sigma_p)\)'s for
\( \sigma_p \in f_{p \, \#}^{-1} (\sigma)\). The representation $\rho_{p,n}$ is a group isomorphism onto its image \(C(\Sigma_p^n, \ZM)\), its inverse is defined as
follows: given \(\varphi = \sum_{\sigma_p \in \Ss_p^n} \varphi_{\sigma_p} \chi_{\sigma_p}\),  where $\varphi_{\sigma_p}$ is an integer,
\(\rho_p^{n \, -1} (\varphi)\) is the group homomorphism from $C_{p,n}$ to $\ZM$ whose value on the basis map $\sigma_p$ is $\varphi_{\sigma_p}$.

\vspace{.1cm}

\noindent Consider the characteristic map $\sigma_p$ of an $n$-simplex $e_p$ in $\Bb_p$. The simplex $e_p$ is contained in some tile $\pi_j$. Viewing $e_p$ as a subset of the
patch $p_j$ in $\RM^d$,  it is possible to define the vector \(x_{\!\sigma_{\!p} \, \partial_i \!\sigma_{\!p}}\), for some $i$ in $1, \cdots n$, that joins the puncture of the
$i$-th face $\partial_i e_p$ (the simplex in $\Bb_p$ whose characteristic map is $\partial_i \sigma_p$) to the puncture of $e_p$.
As a consequence of remark \ref{Ktiling07.rem-orientBp}, those
vectors \(x_{\!\sigma_{\!p} \, \partial_i \!\sigma_{\!p}}\)'s are identical for all $\sigma_p$'s in the preimage of the characteristic map
$\sigma$ of a simplex $e$ on $\Bb_0$, and equal to the vector $x_{\!\sigma \partial_i \!\sigma}$ which defines the operator $\theta_{\!\sigma \partial_i \!\sigma}$ in definition
\ref{Ktiling07.def-ringtrans}.
By analogy, let \(\theta_{\!\sigma_{\!p} \, \partial_i \!\sigma_{\!p}}\) be the operator \(\chi_{\sigma_p} \tra^{x_{\!\sigma \partial_i \!\sigma}} \chi_{\partial_i \sigma_p}\).
With the relation \( \tra^a \chi_\Lambda = \chi_{\tra^a \Lambda} \tra^a\) for
$\Lambda \subset \Omega, a\in \RM^d$,  and lemma \ref{Ktiling07.lem-partitionXip} it is easily seen that $\theta_{\!\sigma \partial_i \!\sigma}$ is the sum
of the \(\theta_{\!\sigma_{\!p} \, \partial_i \!\sigma_{\!p}}\)'s over all \(\sigma_p \in f_{p \,\#}^{-1} (\sigma)\).
Hence PV differential given in equation \eqref{Ktiling07.eq-pimdif} can be written

\begin{equation}
\label{Ktiling07.eq-pimdifBp}
d_{PV}^n = 
\sum_{\sigma_p \in \Ss_{p}^n} \sum_{i=0}^{n} (-1)^i \theta_{\!\sigma_{\!p} \, \partial_i \!\sigma_{\!p}} \,,
\end{equation}

\noindent and is then well defined as a differential from \( C(\Sigma_p^{n-1}, \ZM)\) to \(C(\Sigma_p^n, \ZM) \).

\begin{defini}
\label{Bp}
Let \(C^n_{PV} (p) = C(\Sigma_p^n, \ZM)\), for $n=0, \cdots d$.
The {\rm PV cohomology} of the patch space $\Bb_p$, denoted \(H^\ast_{PV} ( \Bb_0 ; C(\Sigma_p,\ZM) )\), is the homology of the complex
\( \bigl\{ C^n_{PV}(p), d_{PV}^n \bigr\} \).
\end{defini}

\noindent The notation \(H^\ast_{PV} ( \Bb_0 ; C(\Sigma_p,\ZM) )\) requires some comments. The map \(f_p : \Bb_p \rightarrow \Bb_0\) is a
``branched covering'' with discrete ``fibers'', the $\Ss_p^n$'s, that correspond to the ``discrete transversals'' $\Sigma_p^n$'s:

\[
\begin{array}{ccl}
\Sigma_p & \hookrightarrow & \Bb_p \\
 & & \downarrow  \pG_p \\
& & \Bb_0
\end{array}
\]
In analogy with remark \ref{Ktiling07.rem-localcoef}, the PV cohomology
of $\Bb_p$ is analogous to a cohomology of the base space $\Bb_0$ with local coefficients in the $K$-theory of the ``fiber'' $\Sigma_p$.

\begin{proposi}
\label{Ktiling07.prop-pimBp}
The PV cohomology of the patch space $\Bb_p$ is isomorphic to its integer simplicial cohomology: 
\( H^\ast_{PV} ( \Bb_0 ; C(\Sigma_p,\ZM) )  \cong H^\ast (\Bb_p ; \ZM) \).
\end{proposi}

\noindent {\it Proof.} We actually prove a stronger statement: the complexes $C_{PV}(p)^\ast$ and $C_p^\ast$ are chain-equivalent.
As mentioned earlier, $\rho_{p,n}$ in equation \eqref{Ktiling07.eq-reppimcochain} defines a group isomorphism between
the simplicial cochain group $C_p^n$ and the PV cochain group $C_{PV}^n(p)$. Here they are both isomorphic to the direct sum
\( \ZM^{| \Ss_p^n |}\), where $| \Ss_p^n |$ is the cardinality of $\Ss_p^n$, {\it i.e.} the number of $n$ simplices on $\Bb_p$.

\vspace{.1cm}

\noindent  The differential of $\varphi \in C_{PV}^{n-1}(p)$, evaluated on an $n$-simplex $\sigma_p$, reads using \eqref{Ktiling07.eq-pimdifBp}
\(d_{PV}^n \varphi (\sigma_p)=\sum_{\sigma_p \in \Ss_p^n} \bigl(d_{PV}^n \varphi \bigr)_{\sigma_p}  \chi_{\sigma_p} \), with 
\( \bigl(d_{PV}^n \varphi \bigr)_{\sigma_p} = \sum_{i=1}^n (-1)^i \varphi_{\partial_i \sigma_p}\). On the other hand the simplicial differential of 
$\psi \in C^{n-1}_p$ reads \( \delta^n \psi (\sigma_p) = \sum_{i=1}^n (-1)^i \psi(\partial_i \sigma_p)\), and therefore
\( d_{PV}^n \circ \rho_{p,n\!-\!1} = \rho_{p,n} \circ \delta^n , \, n=1, \cdots d\).
Hence the $\rho_{p,n}$'s give a chain map and conjugate the differentials, and thus yield isomorphisms
$\rho_{p,n}^\ast$'s between the $n$-th cohomology groups. \qed

\begin{lemma}
\label{Ktiling07.lem-invlimtrans}
Let \( \bigl\{ \Bb_{l}, f_l \bigr\}_{l\in \NM} \) be a proper sequence of patch spaces of $\Tt$. The following holds: 
\( \Xi_\Delta \cong \varprojlim \bigl( \Ss_l, f_l \bigr)\), and \( C(\Xi_\Delta,\ZM) \cong \varinjlim \bigl( C(\Sigma_l,\ZM) , f^l \bigr)\),
where $f^l$ is the dual map to $f_l$.
\end{lemma}

\noindent {\it Proof.} This is a straightforward consequence of theorem \ref{Ktiling07.thm-invlimBp}. \qed

\begin{proposi}
\label{Ktiling07.prop-pimB0}
Let \( \bigl\{ \Bb_{l}, f_l \bigr\}_{l\in \NM} \) be a proper sequence of patch spaces of $\Tt$. There is an isomorphism:
\[
H_{PV}^{\ast} \bigl( \Bb_0 ; C(\Xi_\Delta,\ZM) \bigr) \cong
\varinjlim \left( H_{PV}^{\ast} \bigl( \Bb_0 ; C(\Sigma_l,\ZM) \bigr) , f_{l}^\ast \right) \,.
\]
\end{proposi}

\noindent {\it Proof.} By the previous lemma \ref{Ktiling07.lem-invlimtrans} the cochain groups $C_{PV}^n$ of the hull
are the direct limits of the cochain groups $C_{PV}^n(l)$ of the patch spaces $\Bb_l$. Let \(f_{l}^{\#}: C_{PV}^n(l) \rightarrow C_{PV}^n (l)\) denote the map induced by $f_l$ on the PV cochain groups.
Since the PV differential $d_{PV}$ 
is the same for the complexes of each patch space $\Bb_l$, it suffices to check that the following diagram
is commutative
\[
\begin{array}{ccccccc}
\cdots & \longrightarrow &  C_{PV}^{n-1}(l) &
\xrightarrow{ \quad d_{P}^{n} \quad } & C_{PV}^{n}(l) &
\longrightarrow & \cdots \\
 & & \downarrow f_{l}^{\#} &  & \downarrow f_{l}^{\#} &  & \\
\cdots & \longrightarrow & C_{PV}^{n-1}(l+1) &
\xrightarrow{ \quad d_{P}^{n} \quad } & C_{PV}^{n}(l+1) &
\longrightarrow & \cdots 
\end{array}
\]
\noindent which is straightforward using the relation \( \bigl( f_l^{\#} \varphi \bigr)_{\partial_j \sigma_l} = \varphi_{ f_{l\, \#}(\partial_j \sigma_l)}\). \qed

\vspace{.1cm}

\noindent {\bf Proof of theorem \ref{Ktiling07.thm-pimcoho}.}
By standard results in algebraic topology \cite{Hatcher02,Spa66}, the \v{C}ech cohomology of the hull is isomorphic to the direct limit of the \v{C}ech cohomologies
of the patch spaces $\Bb_i$,  \( \check{H}^\ast (\Omega ; \ZM) \cong \varinjlim  \left( \check{H}^{\ast} \bigl( \Bb_l ; \ZM \bigr) , f_{l}^\ast \right) \).
By the natural isomorphism between \v{C}ech and simplicial cohomologies for finite $CW$-complexes \cite{Spa66},
\( \check{H}^\ast  \bigl( \Bb_l ; \ZM \bigr) \cong H^\ast \bigl( \Bb_l ; \ZM \bigr)\), and by proposition
\ref{Ktiling07.prop-pimBp}, \( H^{\ast} \bigl( \Bb_l ; \ZM \bigr) \cong H_{PV}^{\ast} \bigl( \Bb_0 ; C(\Sigma_l,\ZM) \bigr) \), so that
\( \check{H}^\ast  \bigl( \Bb_l ; \ZM \bigr) \cong H_{PV}^{\ast} \bigl( \Bb_0 ; C(\Sigma_l,\ZM) \bigr) \). By proposition
\ref{Ktiling07.prop-pimB0}, the direct limit of the $H_{PV}^{\ast} \bigl( \Bb_0 ; C(\Sigma_l,\ZM) \bigr)$'s is the PV cohomology of the hull.
Therefore the integer \v{C}ech cohomology of the hull is isomorphic to the PV cohomology of the hull:
\(\check{H}^\ast(\Omega; \ZM) \cong H^\ast_{PV} (\Bb_0 ;C(\Xi_\Delta,\ZM) )\). \qed

\subsection{Proof of theorem \ref{Ktiling07.thm-SS}}
\label{Ktiling07.ssect-proofSS}

\noindent The proof of theorem \ref{Ktiling07.thm-SS} follows from propositions \ref{Ktiling07.prop-invlimSS} and \ref{Ktiling07.prop-AHSS} below.
First by Thom-Connes isomorphism \cite{Co81},
\( K_{\ast+d} \bigl( C(\Omega) \rtimes \RM^d \bigr) \cong K_{\ast} \bigl( C(\Omega) \bigr) \), and therefore it suffices to build
a spectral sequence that converges to the $K$-theory of the \Cs $C(\Omega)$. This is done by constructing a Schochet spectral sequence
\cite{Schochet81}, $\bigl\{ E^{r,s}_\ast \bigr\}$, associated with an appropriate filtration of $C(\Omega)$. It is shown in proposition
\ref{Ktiling07.prop-invlimSS} that this spectral sequence is the direct limit of spectral sequences, $\bigl\{ E^{r,s}_\ast (l) \bigr\}$,
for the $K$-theory of the \Css $C(\Bb_l)$ of continuous functions on the patch spaces of a proper sequence. Then $\bigl\{ E^{r,s}_\ast (l) \bigr\}$
is shown to be isomorphic to the Atiyah-Hirzebruch spectral sequence \cite{AH61} for the topological $K$-theory of $\Bb_l$ in proposition
\ref{Ktiling07.prop-AHSS}.

\vspace{.1cm}

\noindent For basic definitions, terminology and results on spectral sequences, the reader is referred to \cite{McC}. In section
\ref{Ktiling07.sect-SS} some technical results on exact couples used here are provided.

\vspace{.1cm}

\noindent For $s=0, \cdots d$, let \(\Omega^s=\pG_0^{-1} \bigl( \Bb_0^s \bigr)\) be the lift of the $s$-skeleton of the prototile space, and let
\(I_s = C_0 \bigl( \Omega \setminus \Omega^s \bigr)\). $I_s$ is a closed two-sided ideal of $C(\Omega)$, it consists of functions that vanish on the faces of
dimension $s$ of the boxes of the hull (see section \ref{Ktiling07.ssect-boxhull}). This gives a filtration of $C(\Omega)$:

\begin{equation}
\label{Ktiling07.eq-filtrationhull}
\{0\}=  I_d \hookrightarrow \, I_{d-1} \hookrightarrow \cdots \, I_0 \hookrightarrow I_{-1} \; = \ C(\Omega)  \,.
\end{equation}

\noindent Let \( Q_s = I_{s-1} / I_s \), which is isomorphic to \( C_0 \bigl( \Omega^s \setminus \Omega^{s-1} \bigr) \). Let $K(I)$ and $K(Q)$ be the respective
direct sums of the $K_\epsilon (I_s)$ and $K_\epsilon (Q_s)$ over $\epsilon=0,1$, and $s=0, \cdots d$.
The short exact sequences \( 0 \rightarrow I_s \xrightarrow{\; i_s \; } I_{s-1} \xrightarrow{ \; \pi_s \; } Q_s \rightarrow 0 \), lead,
through long exact sequences in $K$-theory, to the exact couple
\(\PG = \bigl( K(I), K(Q), i, \pi, \partial \bigr)\), where $i$ and $\pi$ are the induced maps and $\partial$ the boundary map in $K$-theory.
Its associated Schochet spectral sequence \cite{Schochet81}, $\bigl\{ E^{r,s}_\ast \bigr\}$, converges to the $K$-theory of $C(\Omega)$:

\begin{equation}
\label{Ktiling07.eq-Schochethull}
\left\{
\begin{array}{ccl}
E^{r,s}_1 & \Rightarrow &  K_{r+s} \bigl( C(\Omega) \bigr) \\
E^{r,s}_1 & = &  K_{r+s} \bigl( Q_s \bigr)
\end{array}
\right. \,.
\end{equation}

\noindent Let  $\Bb_p$ be a patch space  associated with a pattern $\hat{p}$ of $\Tt$. Consider the filtration of the \Cs $C(\Bb_p)$ by the closed
two sided-ideals \(I_s (p)= C_0 \bigl( \Bb_p \setminus \Bb_p^s \bigr)\):

\begin{equation}
\label{Ktiling07.eq-filtrationBp}
\{0\}=  I_d (p) \hookrightarrow \, I_{d-1}(p) \hookrightarrow \cdots \, I_0 (p)  \hookrightarrow I_{-1} (p) \; = \ C(\Bb_p)  \,.
\end{equation}

\noindent Let  \( Q_s(p) = I_{s-1}(p) / I_s (p) \), which is isomorphic to \( C_0 \bigl( \Bb_p^s \setminus \Bb_p^{s-1} \bigr) \). Let $K(I(p))$ and $K(Q(p))$ be the respective
direct sums of the $K_\epsilon (I_s(p))$ and $K_\epsilon (Q_s(p))$ over $\epsilon=0,1$, and $s=0, \cdots d$.
The short sequences \( 0 \rightarrow I_s (p) \xrightarrow{\; i_{p,s}\; } I_{s-1}(p) \xrightarrow{ \; \pi_{p,s} \; } Q_s(p) \rightarrow 0 \), lead to the exact couple
\(\PG(p)= \bigl( K(I(p)), K(Q(p)), i_p, \pi_p, \partial \bigr)\), and its associated Schochet spectral sequence \cite{Schochet81},
$\bigl\{ E^{r,s}_\ast (p) \bigr\}$, converges to the $K$-theory of $C(\Bb_p)$:

\begin{equation}
\label{Ktiling07.eq-SchochetBp}
\left\{
\begin{array}{ccl}
E^{r,s}_1 (p) & \Rightarrow &  K_{r+s} \bigl( C(\Bb_p) \bigr) \\
E^{r,s}_1 (p) & = &  K_{r+s} \bigl( Q_s (p) \bigr)
\end{array}
\right. \,.
\end{equation}

\begin{proposi}
\label{Ktiling07.prop-invlimSS}
Given a proper sequence \( \bigl\{ \Bb_{l}, f_l \bigr\}_{l\in \NM} \) of patch spaces of $\Tt$, the following holds:
\[
 \bigl\{ E^{r,s}_\ast \bigr\} \cong \varinjlim \left( \bigl\{ E^{r,s}_\ast (l) \bigr\}  , f_{l \ast} \right) \,,
\]
\noindent where \(\bigl\{ E^{r,s}_\ast (l) \bigr\}\) is the Schochet spectral sequence \eqref{Ktiling07.eq-SchochetBp}
for $\Bb_l$ corresponding to the patch $p=p_l$.
\end{proposi}

\noindent {\it Proof.} The map \(f_l : \Bb_l \rightarrow \Bb_{l+1}\) induces a morphism of exact couples from \( f_{l\ast} : \PG(l) \rightarrow \PG(l+1)\).
Indeed, consider the following diagram:

\[
\begin{array}{ccccccc}
   K_{\varepsilon} \bigl(I_s(l) \bigr) & \xrightarrow{ \ i_{l,s} \ } &
   K_ {\varepsilon} \bigl( I_{s\!-\!1}(l) \bigr)& \xrightarrow{ \ \pi_{l,s} \ } &
   K _ {\varepsilon} \bigl( Q_s (l) \bigr)&  \xrightarrow{ \ \partial \ } &
   K_{\varepsilon+1} \bigl(I_s (l) \bigr) \\
\downarrow f_{l \ast} && \downarrow  f_{l \ast} &&
 \downarrow f_{l \ast} && \downarrow  f_{l \ast} \\
   K_{\varepsilon} \bigl( I_s (l\!+\!1) \bigr)& \xrightarrow{ \ i_{l\!+\!1,s} \ } &
   K_ {\varepsilon} \bigl(I_{s\!-\!1} (l\!+\!1) \bigr) & \xrightarrow{ \ \pi_{l\!+\!1,s} \ } &
   K _ {\varepsilon} \bigl( Q_s (l\!+\!1) \bigr) &  \xrightarrow{ \ \partial \ } &
   K_{\varepsilon+1} \bigl( I_s (l\!+\!1) \bigr) \\
\end{array}
\]
 
\noindent The left and middle squares are easily seen to be commutative. To check the commutativity of the right square, recall that given a short exact
sequence \( 0 \rightarrow I \rightarrow A \rightarrow A/J \rightarrow 0\) where $A$ is a $C^\ast$-algebra and $J$ a closed two-sided ideal,
the boundary map of an element \( [x] \in K_\varepsilon (A/J) \) is computed via a lift $ z \in A \otimes \Kk $ of $x$.
Let \([x] \in K _ {\varepsilon} \bigl( Q_s (l) \bigr) \). If \(z \in I_{s\!-\!1}(l) \otimes \Kk \) is a lift of $x$, then \( f_{l \ast} z \in I_{s\!-\!1}(l\!+\!1) \otimes \Kk \)
is a lift of $f_{l\ast} x$ and the commutativity of the right square follows.

\vspace{.1cm}

\noindent As a consequence of theorem \ref{Ktiling07.thm-invlimBp}, \( I \cong \varinjlim \{ I(l), f^{\#}_l \}\) and   \( Q \cong \varinjlim \{ Q(l), f^{\#}_l \}\),
therefore \( K(I) \cong \varinjlim \bigl\{ K\bigl( I(l) \bigr), f^{\ast}_l \bigr\} \) and   \( K(Q) \cong \varinjlim \bigl\{ K\bigl(Q(l) \bigr) , f^{\ast}_l \bigr\}\).
Hence \(\bigl\{ \PG(l), f^{\ast}_l \bigr\} \) is a direct system of exact couples and by lemma \ref{Ktiling07.lem-dirlimEC}
\(\PG \cong \varinjlim \bigl\{ \PG(l), f^{\ast}_l \bigr\} \), and the same result on their associated
spectral sequences follows by corollary \ref{Ktiling07.cor-dirlimSS}. \qed

\vspace{.1cm}

\noindent Given a finite $CW$-complex $X$, the Atiyah-Hirzebruch spectral sequence \cite{AH61} for the topological $K$-theory of $X$
is a particular case of Serre spectral sequence \cite{DaKi,McC} for the trivial fibration of $X$ by itself with fiber a point:
\begin{equation}
\label{Ktiling07.eq-AHss}
\left\{
\begin{array}{ccl}
E^{r,s}_{ \! \! \! \! _{\mbox{\tiny \it A\!H}} \, 2} & \Rightarrow &  K^{r+s} ( X ) \\
E^{r,s}_{ \! \! \! \! _{\mbox{\tiny \it A\!H}} \, 2} (p) & \cong &  H^{r} \bigl( X ;  K^s (\cdot) \bigr)
\end{array}
\right. \,,
\end{equation}

\noindent where $\cdot$ denotes a point, so \(K^s (\cdot) \cong \ZM\) if $s$ is even, and $0$ if $s$ is odd.

\vspace{.1cm}

\noindent It is defined on page-$1$ by \(E^{r,s}_{ \! \! \! \! _{\mbox{\tiny \it A\!H}} \, 2} = K^r ( X^s,X^{s-1}) \)  in \cite{AH61}, and then proven that
the page-$2$ is isomorphic to the cellular cohomology of $X$. If $X$ is a locally compact Hausdorff space, the spectral sequence can be rewritten
algebraically, using the isomorphisms \( K^\ast (Y) \cong K_\ast \bigl( C_0 (Y) \bigr) \) and
\(K^\ast (Y,Z) \cong K_\ast \bigl( C_0(Y) / C_0(Z)\bigr)\) for locally compact Hausdorff spaces $Y,Z$.
Consider the Schochet spectral sequence \cite{Schochet81}  for the $K$-theory of the \Cs $C_0(X)$ associated with its filtration by the ideals
$I_s=C_0 (X,X^s)$ of functions vanishing on the $s$-skeleton. Then the spectral sequence associated with the cofiltration of $C_0(X)$ by the ideals
\(F_s = C_0(X) /I_s \cong C_0(X^s)\) turns out to be this algebraic form of Atiyah-Hirzebruch spectral sequence.

\begin{proposi}
\label{Ktiling07.prop-AHSS}
The Schochet spectral sequence \(\bigl\{ E^{r,s}_\ast (p) \bigr\}\) for the $K$-theory of the \Cs $C(\Bb_p)$ is isomorphic to the
Atiyah-Hirzebruch spectral sequence for the topological $K$-theory of $\Bb_p$.
\end{proposi}

\noindent {\it Proof.} By theorem \ref{Ktiling07.thm-filcofil} (section \ref{Ktiling07.sect-SS}) the Schochet spectral sequences built from the
filtration of $A(p) = C(\Bb_p)$ by the ideals $I_s(p)$, namely \( \bigl\{ E^{r,s}_\ast (p) \bigr\} \), and from the cofiltration
of  $A(p)$ by the quotients $F_s (p) = A(p) / I_s(p) \cong C_0 (\Bb_p^s)$, are isomorphic. As remarked above this last spectral sequence
is nothing but the Atiyah-Hirzebruch spectral sequence. \qed

\vspace{.1cm}

\noindent {\bf Proof of theorem \ref{Ktiling07.thm-SS}.} Let \( \bigl\{ \Bb_{l}, f_l \bigr\}_{l\in \NM} \) be a proper sequence of patch spaces of $\Tt$.
By proposition \ref{Ktiling07.prop-AHSS}, the page-$2$ of Schochet spectral sequence \(\bigl\{ E^{r,s}_\ast (l) \bigr\}\) is isomorphic
to the simplicial cohomology of $\Bb_l$: \( E^{r,s}_2 (l) \cong \check{H}^r (\Bb_i ; \ZM )\) for $s$ even and $0$ for $s$ odd. By the natural isomorphism
between simplicial and \v{C}ech cohomologies for $CW$-complexes \cite{Spa66}, it follows that  \( E^{r,s}_2 (l) \cong \check{H}^r (\Bb_i ; \ZM )\) for $s$ even and
$0$ for $s$ odd. By standard results in algebraic topology \cite{Hatcher02,Spa66}
\( \check{H}^\ast (\Omega ; \ZM) \cong \varinjlim  \left( \check{H}^{\ast} \bigl( \Bb_l ; \ZM \bigr) , f_{l}^\ast \right) \), and therefore by proposition
\ref{Ktiling07.prop-invlimSS} the page-$2$ of Schochet spectral sequence for the $K$-theory of $C(\Omega)$ is isomorphic to the integer \v{C}ech
cohomology of the hull: \( E^{r,s}_2 \cong \check{H}^r (\Omega ; \ZM )\) for $s$ even and $0$ for $s$ odd. \qed

\vspace{.3cm}

\section{Appendix: a reminder on spectral sequences}
\label{Ktiling07.sect-SS}

\noindent Some technical results about exact couples are recalled here. 
The formalism allows to state corollary \ref{Ktiling07.cor-extri} which leads to the isomorphism between the Schochet spectral sequence associated with the filtration of a \Cs and the Schochet spectral sequence associated with its corresponding cofiltration in theorem \ref{Ktiling07.thm-filcofil}. The construction of a direct limit of spectral sequences is also recalled.
A few elementary definitions are given to set the terminology. The original reference on exact couples is the work of Massey \cite{Ma52}. The link with
spectral sequences is only mentioned, and the reader is referred to \cite{DaKi,McC} for further details.

\vspace{.1cm}

\noindent An {\em exact couple} is a family $T=(D,E,i,j,k)$ where $D$ and $E$ are abelian groups and $i,j,k$ are
group homomorphisms making the following triangle exact
\begin{equation}
\label{Ktiling07.eq-triangle}
\begin{array}{rcccl}
  D & & \xrightarrow{ \ \ \ i \ \ \ }& & D\\
  k & \nwarrow & & \swarrow & j\\
    & & E & &
  \end{array}
\end{equation}

\noindent In more generality, $D$ and $E$ can be graded modules over a ring, and $i,j,k$ module maps of various degrees.
A morphism $\alpha$ between two exact couples $T$ and $T'$ is a pair of group homomorphisms $(\alpha_D, \alpha_E)$ making the
following diagram commutative
\begin{equation}
\label{Ktiling07.eq-morph}
\begin{array}{lclclcl}
   D &  \xrightarrow{ \ \ \ i \ \ \ }& 
     D &  \xrightarrow{ \ \ \ j \ \ \ } & 
      E &  \xrightarrow{ \ \ \ k \ \ \ } & D \\
\downarrow \alpha_D && \downarrow \alpha_D &&
 \downarrow \alpha_E && \downarrow \alpha_D \\
   D' &  \xrightarrow{ \ \ \ i' \ \ \ } & 
     D' &  \xrightarrow{ \ \ \ j' \ \ \ } & 
      E' & \xrightarrow{ \ \ \ k' \ \ \ }& D'
  \end{array}
\end{equation}

\noindent The composition map $d=j\circ k : E\mapsto E$ is called the {\it differential} of the exact couple. 
Since \eqref{Ktiling07.eq-triangle} is exact, it follows that $d^2= j\circ (k\circ j)\circ k =0$. Let then $H_d(E)$ be the homology of the 
complex $E\stackrel{d}{\longrightarrow}E$. Then the following theorem holds.

\begin{theo}
\label{Ktiling07.th-dertria}
\begin{enumerate}[(i)]

\item There is a {\em derived exact couple}
\begin{equation}
\label{exactcouples06.eq-dertria}
\begin{array}{rcccl}
  D^{(1)}=i(D) & &\xrightarrow{ \quad \quad \quad  i^{(1)} \quad \quad  \quad }&
                 & D^{(1)}=i(D)\\
   k^{(1)} &  \nwarrow & & \swarrow  &j^{(1)}\\
           &  &E^{(1)}=H_d(E)&    & 
  \end{array}
\end{equation}
\noindent defined by $i^{(1)} =i \vert_{i(D)}$, 
$j^{(1)}(i(x))=j(x) + \Ima(d)$ and $k^{(1)}(e+\Ima(d))=k(e)$. 

\item The derivation \( T \rightarrow T^{(1)}\) is a functor on the category of exact couples (with morphisms of exact couples).
\end{enumerate}
\end{theo}

\noindent The $n$-th iterated derived couple of $T$ is denoted \(T^{(n)} = (D^{(n)}, E^{(n)}, i^{(n)}, j^{(n)}, k^{(n)})\). 

\begin{rem}
\label{Ktiling07.rem-dermorph}
{\em Point (ii) in theorem \ref{exactcouples06.eq-dertria} implies that a morphism of exact couples  \( \alpha :  T \rightarrow T' \) induces derived morphisms  $\alpha^{(n)}$ between the derived
couples: \( \alpha_{D}^{(n)} = \alpha_{D}^{(n-1)} \vert_{ D^{(n)} } \) is the restriction of $\alpha_{D}^{(n-1)}$ to $D^{(n)}=i^{(n-1)}( D^{(n-1)} )$, 
and \( \alpha_{E}^{(n)} = ( \alpha_{E}^{(n-1)} )_{\ast} \) is the induced map of $\alpha_{E}^{(n-1)}$  on homology.
Also $\alpha_E$ conjugates the differentials, namely \( \alpha_E^{(n)} \circ d^{(n)} = d'^{(n)} \circ \alpha_E^{(n)} \) for all $n$.}
\end{rem}

\noindent The {\it spectral sequences} considered here come from exact couples. Given an exact couple
$(E,D,i,j,k)$, its associated spectral sequence is the family $( E_l, d_l)_{l\in\NM}$ where
$E_l =E^{(l)}$ is the derived $E$-term, and called the {\it $E_l$-page} or simply {\it page-$l$} of the spectral sequence,
and $d_l = d^{(l)}$ is its differential. As noted above those $E$-terms are generally
(bi)graded modules. A morphism of spectral sequence \(\beta : (E_l, d_l)_{l\in \NM} \rightarrow (E'_l, d'_l)_{l\in \NM}\) is given by the
$\alpha_E$-morphism of the exact couples, and can be seen as a sequence of {\it chain maps} 
\(\beta_l : E_l \rightarrow E'_l\) for each $l$, {\it i.e.} commute with the differentials: \(d_l \beta_l = \beta_l d'_l\).

\begin{defini}
\label{Ktiling07.def-varEC}
\begin{enumerate}[(i)]
\item An exact couple $T=(D,E,i,j,k)$ is said to be {\rm trivial} if $E=0$.
\item Two exact couples $T$ and $T'$ are said to be {\rm equivalent} if there is a morphism $\alpha : T \rightarrow T'$
such that $\alpha_{E}$ is an isomorphism.

\end{enumerate}
\end{defini}

\begin{rem}
\label{Ktiling07.rem-trivEC}
{\em 
\begin{enumerate}[(i)]

\item If an exact couple is trivial, then it is of the form $(D,0,{\rm id}, 0, 0)$. 

\item A trivial exact couple is identical to its derived couple.

\item An equivalence between two exact couples is equivalent to an isomorphism between their associated spectral sequences.

\end{enumerate}}
\end{rem}

\begin{defini}
\label{Ktiling07.def-cvEC}
\begin{enumerate}[(i)]

\item An exact couple is said to converge whenever there is an $L\in \NM$ such that for $l\geq L$ the 
$l$-th derived couples are trivial. 

\item A spectral sequence is said to converge if its associated exact couple converges.
\end{enumerate}
\end{defini}

\noindent Given a morphism of exact couples \( \alpha : T \rightarrow T'\), its kernel \( \Ker \alpha = ( \Ker \alpha_D, \Ker \alpha_E, i,j,k)\) and image
\( {\rm Im} \alpha = ( {\rm Im} \alpha_D, {\rm Im} \alpha_D, i', j', k')\) are $3$-periodic complexes (by commutativity of diagram \eqref{Ktiling07.eq-morph}),
but no longer exact couples in general.

\begin{lemma}
\label{Ktiling07.prop-ESEC}
Given an exact sequence of exact couples
\[
\cdots \rightarrow 
T_{m-1} \xrightarrow{\alpha_{m-1}} 
T_{m} \xrightarrow{\alpha_{m}} 
T_{m+1} \xrightarrow{\alpha_{m+1}} 
T_{m+2} \rightarrow  
\cdots \,,
\]
that is  \( {\rm Im} \alpha_{m-1} = \Ker  \alpha_{m} \) for all $m$,  the following holds:
if $T_{m-1}$ and $T_{m+2}$ are trivial, then $T_{m}$ and $T_{m+1}$ are equivalent.
\end{lemma}

\noindent {\bf Proof.}  Using Remark \ref{Ktiling07.rem-trivEC} the couples can be  rewritten \(T_{m-1} =(D_{m-1}, 0, {\rm id}, 0,0) \) and 
\( T_{m+2}=  (D_{m+2}, 0, {\rm id}, 0,0) \), and therefore \( \alpha_{m-1}=(\alpha_{D_{m-1}},0)\), and \(\alpha_{m+1} = (\alpha_{D_{m+1}}, 0)\).
The exact sequence for the $E$ terms then reads 
\( \cdots \rightarrow 0 \rightarrow E_m \xrightarrow{ \alpha_{E_m} } E_{m+1} \rightarrow 0 \rightarrow \cdots\), and 
therefore $\alpha_{E_m}$ is an isomorphism and gives an equivalence between $T_m$ and $T_{m+1}$. \qed
\begin{coro}
\label{Ktiling07.cor-extri}
Given an exact triangle of exact couples, one of which is trivial, then the two others are equivalent. 
\end{coro}

\begin{defini}
\label{Ktiling07.def-dirlimEC}
A direct system of exact couples  \(\bigl\{ T_l, \alpha_{lm} \bigr\}_I\), is given by an directed set $I$, and a family of exact couples
$T_l$'s and morphisms of exact couples  \( \alpha_{lm} : T_l \rightarrow T_m\) for $l\le m$ (with $\alpha_{ll}$ the identity), such that given
$l,m \in I$ there exists $n\in I$, $n\ge l,m$ with \( \alpha_{ln} = \alpha_{mn} \circ \alpha_{lm}\).
\end{defini}

\noindent If $T_l$ is written \( (D_l, E_l, i_l,j_l,k_l)\), the definition of a direct system of associated spectral sequences is given similarly
by keeping only the data of the $E_l$-terms, their differentials $j_l \circ k_l$ and the morphisms $\alpha_{E_l}$.

\vspace{.1cm}

\noindent Let $R$ be a commutative ring. Recall that the direct limit of a directed system of $R$-modules, \(\bigl\{ M_l, \sigma_{lm} \bigr\}_I\) is given
by the quotient of the direct product \( \bigoplus_{I} M_l \) by the $R$-module generated by all elements of the form \(a_l - \sigma_{lm}(a_l)\)
for \(a_l \in M_l\), where each $R$-module $M_l$ is viewed as a submodule of \( \bigoplus_{I} A_l \).

\begin{lemma}
\label{Ktiling07.lem-dirlimEC}
Let \(\bigl\{ T_l, \alpha_{lm} \bigr\}_I\) be a direct system of exact couples, and write the couples as \( T_l = (D_l, E_l, i_l,j_l,k_l)\), and the morphisms as
\( \alpha_{lm} = (\alpha_{D_{lm}}, \alpha_{E_{lm}})\).
There is a direct limit exact couple \(T = \varinjlim \bigl\{ T_l, \alpha_{lm} \bigr\}\), given by \(T = (D,E,i,j,l)\) with
\( D = \varinjlim \bigl\{ D_l, \alpha_{D_lm} \bigr\}\), \( E = \varinjlim \bigl\{ E_l, \alpha_{E_lm} \bigr\}\), and
$i,j,k$ the maps induced by the $i_l,j_l,k_l$'s.
\end{lemma}

\noindent  {\it Proof.} The direct limits $D$ and $E$ are well-defined and it suffices to give the expressions of the module homomorphisms $i,j,k$
and show that $T$ is exact. Let $d$ be in $D$, {\it i.e.} it is the class $[d_l]$ for some $d_l \in D_l$, then $i(d) = [i_l (d_l)]$. If $d=[d_m]$ for some $m\ge l$
then \(d_m = \alpha_{D_{lm}}(d_l)\) and \( [i_m (d_m)] = [ i_m \circ \alpha_{D_{lm}} (d_l)] = [\alpha_{D_{lm}} \circ i_l(d_l)]\) because $\alpha_{D_{lm}}$
is a morphism of exact couple (commutativity of diagram \eqref{Ktiling07.eq-morph}), and therefore \( [i_m(d_m) = [i_l(d_l)]\) and $i$ is well-defined.
Similarly for $d=[d_l]$ in $D$, \(j(d) = [j_l(d_l)]\) in $E$, and for $e=[e_l]$ in $E$, \(k(e)=[k_l(e_l)]\) in $D$, and are well-defined. 

\noindent This proves also that
$T$ is a $3$-periodic complex since the compositions $j\circ i$, $k \circ j$ and $i \circ k$ involves the compositions of the $i_l, j_l, k_l$, and are thus zero.
Let $d \in \Ker j$ and $d=[d_l]$, then \(j(d)= [j_l(d_l)]=0\) and by exactness of $T_l$ there exists $d'_l \in D_l$ such that $d_l=i_l(d'_l)$, let then
$d'=[d'_l]$ in $D$ to have $d=[i_l(d'_l)] = i(d')$ and therefore $\Ker j = {\rm Im} i$. The other two relations \( \Ker i = {\rm Im} k\) and \( \Ker k = {\rm Im} j\) are
proven similarly, and this shows that $T$ is exact.
\qed

\begin{coro} 
\label{Ktiling07.cor-dirlimSS}
The result of lemma \ref{Ktiling07.lem-dirlimEC} for direct limit of exact couples also holds for direct limit of associated spectral sequences.
\end{coro}

\vspace{.1cm}

\noindent Let $A$ be a \Cs, and assume there is a finite filtration by closed two-sided ideals:
\begin{equation}
\label{Ktiling07.eq-fil}
 \{0\}=  I_d \stackrel{i_d}{\hookrightarrow} \, I_{d-1} \stackrel{i_{d\!-\!1}}{\hookrightarrow} \cdots \, I_0 \stackrel{i_0}{\hookrightarrow} I_{-1} \; = \ A \ \,.
\end{equation}
\noindent There is an associated cofiltration of $A$ by the quotient \Css $F_p = A / I_p$:
\begin{equation}
\label{Ktiling07.eq-cofil}
A \ = F_d \stackrel{\rho_d}{\twoheadrightarrow} F_{d-1} \stackrel{\rho_{d\!-\!1}}{\twoheadrightarrow} \cdots F_0 \stackrel{\rho_0}{\twoheadrightarrow} F_{-1} = \{ 0 \} \ \,.
\end{equation}

\noindent Let $Q_p=I_{p-1} / I_p$ be the quotient \CS. There are short exact sequences:
\begin{subequations}
\label{Ktiling07.eq-SESgrp}
\begin{align}
0  \longrightarrow \; I_p \, \xrightarrow{ \ \ i_p \ \ }  I_{p-1}  \xrightarrow{ \ \ \pi_p \ \ } \ Q_p \ \longrightarrow  0 
\label{Ktiling07.eq-SESfiltration}\\
0  \longrightarrow  Q_p \xrightarrow{ \ \ j_p \ \ } \  F_p \  \xrightarrow{ \ \ \rho_p \ \ }  F_{p-1}  \longrightarrow  0 
\label{Ktiling07.eq-SEScofiltration} \\
0  \longrightarrow \; I_p \, \xrightarrow{ \ \ l_p \ \ }  \ \, A  \ \, \xrightarrow{ \ \ \sigma_p \ \ } \  F_p \, \ \longrightarrow  0 
\label{Ktiling07.eq-SESlink}
\end{align}
\end{subequations}

\noindent In \eqref{Ktiling07.eq-SESfiltration} $i_p$ is the canonical inclusion, and $\pi_p$ the quotient map \( \pi_p (x) = x + I_p\). 
In \eqref{Ktiling07.eq-SESlink} $l_p = i_p \circ i_{p-1} \circ \cdots i_{0}$
is the canonical inclusion, and $\sigma_p$ the quotient map \(\sigma_{p} ( x ) = x + I_p\).
And in \eqref{Ktiling07.eq-SEScofiltration} $j_p$ is the canonical inclusion \( j_p ( x + I_p ) = l_{p-1}(x) + I_p \), and $\rho_p$ the 
quotient map \( \rho_p ( x + I_p ) = x + I_{p-1}\).

\vspace{.1cm}

\noindent Associated with the short exact sequences of \Css \eqref{Ktiling07.eq-SESfiltration}  and \eqref{Ktiling07.eq-SEScofiltration} there are a long (6-term periodic) exact sequences in $K$-theory that can be written in exact couples: 

\begin{subequations}
\begin{align}
T_I \ : \ 
\begin{array}{rcccl}
K(I) && \xrightarrow{ \ \ \ \ i \ \  \ \ }  && K(I) \\
\partial &\nwarrow &  & \swarrow & \pi \\
& & K(Q) & &
\end{array} 
\quad \quad \text{with} \quad K(I) = \bigoplus_{p=-1}^d \bigoplus_{\varepsilon=0,1} K_\varepsilon (I_p) 
\label{Ktiling07.eq-ECfiltration}\\
T_F \ : \ 
\begin{array}{rcccl}
\!  K(F) \!  && \xrightarrow{ \ \ \ \ \rho \ \  \ \ }  && \!  \! K(F) \!  \\
j &\nwarrow &  & \swarrow & \partial \\
& & K(Q) & &
\end{array}
\quad \quad \text{with} \quad K(F) = \bigoplus_{p=-1}^d \bigoplus_{\varepsilon=0,1} K_\varepsilon (F_p) 
\label{Ktiling07.eq-ECcofiltration}
\end{align}
\end{subequations}

\noindent with $K(Q) = {\displaystyle \bigoplus_{p=-1}^d \bigoplus_{\varepsilon=0,1}} K_\varepsilon (Q_p)$, and
$i$, $\pi$, $j$ and $\rho$ the induced maps on K-theory.

\vspace{.1cm}

\noindent Since the filtration \eqref{Ktiling07.eq-fil} and the cofiltration \eqref{Ktiling07.eq-cofil} are finite, both exact couple converge to the $K$-theory of $A$.
\begin{theo}
\label{Ktiling07.thm-filcofil}
The spectral sequence for the $K$-theory of a \Cs $A$, associated with a filtration of $A$ by ideals $I_p$ as in \eqref{Ktiling07.eq-fil}, and the spectral sequence associated with its
corresponding cofiltration by quotients $F_p = A /I_p$ as in \eqref{Ktiling07.eq-cofil}, are isomorphic.
\end{theo}

\noindent {\it Proof.} By remark \ref{Ktiling07.rem-trivEC} (iii) it is sufficient to prove that the exact couples $T_I$ \eqref{Ktiling07.eq-ECfiltration} and 
$T_F$ \eqref{Ktiling07.eq-ECcofiltration} are equivalent. 
By corollary \ref{Ktiling07.cor-extri} it is also sufficient to prove that $T_I$ and $T_F$ fit into an exact triangle with a trivial exact couple.

\noindent Let $T_A$ be the the trivial exact couple \( T_A = (K(A),0, {\rm id}, 0, 0 ) \). We prove that the short exact sequences 
\eqref{Ktiling07.eq-SESlink} induce the following exact triangle of exact couples
\begin{equation}
\label{Ktiling07.eq-extri}
\begin{array}{rcccl}
  T_I  && \xrightarrow{ \ \ \ \ l \ \  \ \ }  &&  T_A  \\
\partial &\nwarrow &  & \swarrow & \sigma \\
& & T_F & &
\end{array}
\end{equation}

\noindent The exactness of the triangle comes from the exactness of the K-theory functor, and one has now to verify that the applications $l$, $\sigma$, and $\partial$ define morphisms of exact couples. 

\vspace{.1cm}

\noindent Those proofs are similar. We give the details for \( (l,0): T_I \rightarrow T_A \), where the application $0$ is the induced quotient map of the trivial short exact sequence \( 0 \longrightarrow Q_p \xrightarrow{ \  \ {\rm id}  \ \ }  
Q_p \xrightarrow{ \ \ 0 \ \  } 0 \longrightarrow 0 \), {\it i.e.} we prove that the following diagram is commutative ($\varepsilon$ is $0$ or $1$):

\begin{equation}
\label{Ktiling07.eq-morph1}
\begin{array}{ccccccc}
   K_{\varepsilon} (I_p) & \xrightarrow{ \ \  i_{p \ast}  \ \ } & 
   K_ {\varepsilon} (I_{p-1})& \xrightarrow{ \ \ \pi_{p\ast}  \ \ } & 
   K _ {\varepsilon} (Q_p)&  \xrightarrow{ \ \ \ \partial \ \ \ } & 
   K_{\varepsilon+1} (I_p) \\
\downarrow l_{p\ast} && \downarrow l_{p-1 \ast} &&
 \downarrow 0 && \downarrow  l_{p\ast} \\
   K_{\varepsilon} (A) & \xrightarrow{ \ \ \ {\rm id}_\ast \ \ \ } & 
   K_ {\varepsilon} (A)& \xrightarrow{ \ \ \ 0 \ \ \ } & 
   0 &  \xrightarrow{ \ \ \ \partial \ \ \ } & 
   K_{\varepsilon+1} (A) \\   
\end{array}
\end{equation}

\noindent The middle square is commutative since the maps lead to $0$. The commutativity of the left square comes from the functoriality of 
K-theory: since \( l_p = l_{p-1} \circ i_p \), it follows that \( {\rm id}_\ast l_{p\ast} = l_{p\ast} = l_{p-1\ast} i_{p\ast} \). For the right square \( l_{p\ast} = {\rm id}_\ast l_{p \ast} = l_{p-1\ast} i_{p\ast} \) by
commutativity of the left square, and thus \( l_{p\ast} \partial = l_{p-1\ast} i_{p\ast} \partial = 0\)
because \( i_{p\ast} \partial = 0\) by exactness of the long exact sequence in K-theory induced by \eqref{Ktiling07.eq-SESfiltration}. \qed



\begin{thebibliography}{99}
\bibliographystyle{unsrt}

\bibitem{AP98} Anderson J. E., Putnam I. F., ``Topological invariants for substitution tilings and their associated $C^\ast$-algebra'', in {\it Ergod. Th, \& Dynam. Sys.}, {\bf 18} (1998), 509-537.

\bibitem{AH61} Atiyah M.F., Hirzebruch F., ``Vector bundles and homogeneous spaces'', {\it Proc. Symp. Pure Math.} {\bf 3} 7-38, Providence: Amer. Math.  Soc. (1961).

\bibitem{AS63} Atiyah M.F., Singer I.M., ``The index of elliptic operators on compact manifolds", {\it Bull. Amer. Math. Soc.} {\bf 69} (1963), 422-433.

\bibitem{AS68} Atiyah M.F., Segal G.B., Singer I.M. ``The inder of elliptic operators. I, II, III", {\it Ann. of Math. (2)} {\bf 87} (1968) 484-530, 531-445, 546-604.

\bibitem{Be82} Bellissard J., ``Schr\"{o}dinger's operators with an almost periodic potential : an over\-view'', in {\it Lecture Notes in Phys.}, {\bf 153}, Springer Verlag, Berlin Heidelberg, New-York, (1982).

\bibitem{Be86} Bellissard J., ``$K$-Theory of $C^{\ast}$-algebras in Solid State Physics'', in {\it Statistical Mechanics and Field Theory,  Mathematical Aspects}, T.C. Dorlas, M.N. Hugenholtz \& M. Winnink, {\it Lecture Notes in Physics}, {\bf 257} (1986), 99-156.

\bibitem{BIST89} Bellissard J., Iochum B., Scoppola E., Testard D., ``Spectral Properties of One Dimensional Quasi-Crystals'', {\it Commun. Math. Phys.} {\bf 125} (1989), 527-543.

\bibitem{Be91} Bellissard J., Bovier A., Ghez, J.M., ``Gap Labelling Theorem for one Dimensional Schr\"odinger Operators'', {\it Rev. Math. Phys.}, {\bf 4} (1992), no. 1, 1-37.

\bibitem{Be93} Bellissard J., ``Gap Labelling Theorems for Schr\"{o}dinger's Operators'', in {\it From Number Theory to Physics}, pp.538-630, Les Houches March 89, Springer, J.M. Luck, P. Moussa \& M. Waldschmidt Eds., (1993).

\bibitem{BCL98} Bellissard J., Contensou E., Legrand A., ``$K$-th\'eorie des quasi-cristaux, image par la trace: le cas du r\'eseau octogonal'', {\it C. R. Acad. Sci. Paris S\'er. I Math.}, {\bf 326} (1998), 197-200.

\bibitem{BHZ00} Bellissard J., Hermmann D., Zarrouati M., ``Hull of Aperiodic Solids and Gap labeling Theorems'', in {\it Directions in Mathematical Quasicrystals}, CRM Monograph Series, {bf 13}, 207-259, M.B. Baake \& R.V. Moody Eds., AMS Providence, (2000).

\bibitem{BKL01} Bellissard J., Kellendonk J., Legrand A., ``Gap-labelling for three-dimensional aperiodic solids'', {\it C. R. Acad. Sci. Paris, S\'er. I Math.}, {\bf 332} (2001), 521-525.

\bibitem{BBG06} Bellissard J., Benedetti R., Gambaudo J.-M., ``Spaces of Tilings, Finite Telescopic Approximations and Gap-labelling'', {\it Commun. Math. Phys.}, {\bf 261} (2006), 1-41.

\bibitem{BG03} Benedetti R., Gambaudo J.-M., ``On the dynamics of $\mathbb G$-solenoids. Applications to Delone sets'', {\it Ergodic Theory Dynam. Systems}, {\bf 23} (2003), 673-691.

\bibitem{Bla98} Blackadar B. {\it $K$-Theory for Operator Algebras}, 2nd Ed., Cambridge University Press, Cambridge (1998).

\bibitem{BS50} Borel A., Serre J.-P., ``Impossibilit\'e de fibrer un espace euclidien par des fibres compactes, {\it C. R. Acad. Sci. Paris}, {\bf 230} (1950), 2258-2260.

\bibitem{Br72} Bratelli O., ``Inductive limits of finite dimensional $C^\ast$-algebras'', {\it Trans. Amer. Math. Soc.}, {\bf 171} (1972), 195-234.

\bibitem{CW78} Claro~F.~H., Wannier~G.~H., ``Closure of Bands for Bloch Electrons in a Magnetic Field, {\it Phys. Stat. Sol. B}, {\bf 88} (1978), K147-K151.

\bibitem{CMS73} Coburn L.A., Moyer R.D., Singer I.M., ``\Cs of almost periodic pseudo-differential operators'', {\it Acta Math.},{\bf 130} (1973), 279-307.

\bibitem{Co79} Connes A., ``Sur la th\'eorie non commutative de l'int\'egration", {\it Alg\`ebres d'op\'erateurs. (S\'em. Les Plans-sur-Bex)} (1978), 19-143.

\bibitem{Co81} Connes A., ``An analogue of the Thom isomorphism for crossed products of a $C^\ast$-algebra by an action of R'', {\it Adv. in Math.} {\bf 39} (1981), no. 1, 31-55.

\bibitem{Co82} Connes, A., ``A survey of foliations and operator algebras.  Operator algebras and applications, Part I'', (Kingston, Ont., 1980),  pp. 521-628, {\it Proc. Sympos. Pure Math.}, {\bf 38}, Amer. Math. Soc., Providence, R.I., 1982.

\bibitem{Co90} Connes A., {\it G\'eom\'etrie non commutative}, InterEditions, Paris, (1990).

\bibitem{Co94} Connes A., {\it Noncommutative Geometry}, Academic Press, San Diego, (1994).

\bibitem{DaKi} Davis J.F., Kirk P., {\it Lecture Notes in Algebraic Topology}, AMS Graduate Studies in Mathematics {\bf 35} (2001)

\bibitem{vEl94} Van Elst A., ``Gap-labelling theorems for Schr\"odinger operators on the square and cubic lattice''. {\it Rev. Math. Phys.}, {\bf 6} (1994), 319-342.

\bibitem{FH99} Forrest A. H., Hunton J. R., ``The cohomology and $K$-theory of commuting homeomorphisms of the Cantor set'', {\it Ergodic Theory Dynam. Systems}, {\bf 19} (1999), 611-625.

\bibitem{HFK02} Forrest A. H., Hunton J. R., Kellendonk J., ``Topological Invariants for Projection Method Patterns'', in {\it Memoirs of the American Mathematical Society}, {\bf 159} (2002),  no. 758,

\bibitem{GHK05} G\"ahler F., Hunton J., Kellendonk J., ``Torsion in Tiling Homology and Cohomology'', at arXiv.com {\tt math-phys/0505048}, (May 2005).

\bibitem{GM06} Gambaudo J.-M., Martens M., ``Algebraic topology for minimal Cantor sets''. {\it Ann. Henri Poincar\'e}, {\bf 7} (2006), 423-446.

\bibitem{Gh99} Ghys~E., ``Laminations par surfaces de Riemann'', in {\sl Dynamique et g\'eom\'etrie complexes}, {\it Panoramas \& Synth\`eses}, {\bf 8} (1999), 49-95.

\bibitem{Gro58} Grothendieck A., ``La th\'eorie des classes de Chern'', {\it Bull. Soc. Math. France}, {\bf 86} (1958), 137-154.

\bibitem{GruShep87} B. Gr\"unbaum, G.C. Shephard, {\it Tilings and Patterns}, W.H. Freemand and Co, New York, 1st Ed. (1987).

\bibitem{Hatcher02} Hatcher, A., {\it Algebraic Topology},  1st ed., Cambridge University Press, (2002), ( {\em and available online at} http://www.math.cornell.edu/~hatcher/\#ATI

\bibitem{Hi54} Hirzebruch F., ``Arithmetic genera and the theorem of Riemann-Roch for algebraic varieties'', {\it Proc. Nat. Acad. Sci. U. S. A.}, {\bf 40} (1954), 110-114.

\bibitem{Hi56} Hirzebruch F., {\it Neue topologische Methoden in der algebraischen Geometrie}, Ergebnisse der Mathematik und ihrer Grenzgebiete (N.F.), Heft 9. Springer-Verlag, Berlin-Gšttingen-Heidelberg, 1956. {\it Topological Methods in Algebraic Geometry}, Reprint of the third Edition, Springer-Verlag, Berlin-Heidelberg-New-York, 1978.

\bibitem{Hof76} Hofstadter~D.~R., ``Energy levels and wave functions of Bloch electrons in a rational or irrational magnetic field'', {\it Phys. Rev. B}, {\bf 14} (1976) 2239-2249.

\bibitem{JM82} Johnson R., Moser, J., ``The rotation number for almost periodic potentials''. {\it Comm. Math. Phys.}, {\bf 84} (1982), 403-438.

\bibitem{Ka80} Kasparov, G. G., ``The operator $K$-functor and extensions of $C\sp{*} $-algebras'', {\it Izv. Akad. Nauk SSSR}, {\bf 44} (1980), 571-636.

\bibitem{Ka88} Kasparov G.G., ``Equivariant $KK$-theory and the Novikov conjecture'', {\it Inv. Math.} {\bf 91} (1988), 147-201.

\bibitem{Ke95} Kellendonk J., ``Noncommutative geometry of tilings and gap labelling''. {\it Rev. Math. Phys.}, {\bf 7} (1995), 1133-1180.

\bibitem{Ke03} Kellendonk J., ``Pattern-equivariant functions and cohomology'', {\it J. Phys. A}, {\bf 36} (2003), 5765-5772.

\bibitem{KP06} Kellendonk J., Putnam I.~F., ``The Ruelle-Sullivan map for actions of $\RM^n$'', {\it Math. Ann.}, {\bf 334} (2006), 693-711.

\bibitem{Ko47} Koszul J.-L., ``Sur les op\'erateurs de d\'erivation dans un anneau'', {\it C. R. Acad. Sci., Paris}, {\bf 225} (1947), 217-219.

\bibitem{La99A} Lagarias J.~C., ``Geometric models for quasicrystals I. Delone sets of finite type''. {\it Discrete Comput. Geom.}, {\bf 21} (1999), 161-191.

\bibitem{La99B} Lagarias J.~C., ``Geometric models for quasicrystals. II. Local rules under isometries''. {\it Discrete Comput. Geom.}, {\bf 21} (1999), 345-372.

\bibitem{LP03} Lagarias J.~C., Pleasants P.~A.~B., ``Repetitive Delone sets and quasicrystals''. {\it Ergodic Theory Dynam. Systems}, {\bf 23} (2003), 831-867.

\bibitem{Le46} Leray J., ``Structure de l'anneau d'homologie d'une repr\'esentation'', {\it C. R. Acad. Sci., Paris}, {\bf 222} (1946), 1419-1422.

\bibitem{Ma52} Massey W.~S., ``Exact Couples in Algebraic Topology (Parts I \& II)'', {\it Ann. Math.}, {\bf 56}, (1952), 363-396; ``Exact Couples in Algebraic Topology (Parts III, IV \& V)'', {\it Ann. Math.}, {\bf 57} (1953), 248-286.

\bibitem{McC} McCleary J., {\it A User's Guide to Spectral Sequences}, Cambridge Studies in Advanced Mathematics {\bf 58}, 2nd Ed., Cambridge University Press, Cambridge (2001).

\bibitem{Mey72} Meyer Y., {\it Algebraic Number and Harmonic Analysis}, North-Holland, Amsterdam (1972).

\bibitem{MS88} Moore~C.~C., Schochet~C., {\sl  Global Analysis on Foliated Spaces}, Math. Sci. Res.  Inst. Publ. {\bf No. 9}, New York; Springer-Verlag,  (1988).

\bibitem{Mos81} Moser J., ``An example of a Schr\"odinger equation with almost periodic potential and nowhere dense spectrum'', {\it Comment. Math. Helv.}, {\bf 56} (1981), 198-224.

\bibitem{Pa65} Palais R., {\it Seminar on the Atiyah-Singer index theorem}, Ann. of Math. Studies, no 57, Princeton Univ. Press, Princeton, N.J., 1965

\bibitem{PV80} Pimsner M., Voiculescu D., ``Exact sequences for $K$-groups and Ext groups of certain cross-product $C^*$-algebras''. {\it J. Operator Theory}. {\bf 4} (1980), 93-118.

\bibitem{Pi83} Pimsner M.~V., ``Ranges of traces on $K_0$ of reduced crossed products by free groups'', in {\it Operator algebras and their connections with topology and ergodic theory (Bu\c steni, 1983)}, 374-408, {\it Lecture Notes in Math.}, {\bf 1132}, Springer, Berlin, (1985).

\bibitem{Pen74} Penrose R., ``The role of aesthetics in pure and applied mathematical research'' {\it Bull. Inst. Math. Appl.} {\bf 10} (1974), 55-65.

\bibitem{QU87} Queffelec, M., {\it Substitution dynamical systems-Spectral analysis}. {\it Lecture Notes in Math. Springer}, Berlin Heidelberg New York (1987), vol. {\bf 1294}.

\bibitem{Rad94} Radin C., ``The pinwheel tiling of the plane'', {\it Ann. of Math. (2)} {\bf 139} (1994), 257-264.

\bibitem{Rad99} Radin C., {\it Miles of Tiles}, AMS, Student Mathematical Library, vol. 1, (1999).

\bibitem{Re80} Renault J., {\it A groupoid approach to \Css} Lecture Notes in Math. {\bf 793} Berlin-Heidelberg-New-York: Springer-Verlag (1980).

\bibitem{Ri82} Rieffel M.A., ``Morita equivalence for operator algebras'', {\it Operator algebras and applications}, {\it Proc. Symposia Pure Math.} {\bf 32} (1982), Part 1, 285-298.

\bibitem{Sa03} Sadun L., ``Tiling spaces are inverse limits'', {\it J. Math. Phys.}, {\bf 44} (2003), 5410-5414.

\bibitem{SW03} Sadun L., Williams R.~F., ``Tiling spaces are Cantor set fiber bundles'', {\it Ergodic Theory Dynam. Systems}, {\bf 23} (2003), 307-316.

\bibitem{Sa05} Sadun L., ``Explicit computation of the cohomology for the chair tiling'', talk delivered at Banff Research Station, July 2005.

\bibitem{Sa06} Sadun L., ``Pattern-Equivariant Cohomology with Integer Coefficients'', at arXiv.com {\tt math.DS/0602066}, Feb. 2006. Available at\\
ftp://ftp.ma.utexas.edu/pub/papers/sadun/2006/integer4.pdf

\bibitem{Se51} Serre J.-P., ``Homologie singuli\`ere des espaces fibr\'es. Applications'', {\it Ann. of Math.}, {\bf 54} (1951), 425-505.

\bibitem{SBGC84} Shechtman, D., Blech, I., Gratias, D., Cahn, J.V., ``Metallic phase with long range orientational order and no translational symmetry''. {\it Phys. Rev. Lett.}, {\bf 53} (1984), 1951-1953.

\bibitem{Schochet81} Schochet, C., ``Topological methods for C*-algebra I: Spectral sequences'', {\it Pacif. J. Math.} {\bf 96} (1981), no. 1, 193-211.

\bibitem{Si70} Singer I.M., ``Future extensions of index theory and elliptic operators", {\it Proc. Sympos., Princeton Univ., Princeton, N.J.,}(1970) 171-185 

\bibitem{Spa66} Spanier E., {\it Algebraic Topology}, McGraw-Hill, 1966 (reprinted by Springer-Verlag).

\bibitem{Su79} Subin M.A., ``Spectral theory and the index of elliptic operators with almost-periodic coefficients" (Russian) {\it Uspeki Mat. Nauk} {\bf 34} (1979), no. 2(206), 95-135.

\bibitem{Sut87} S\"uto, A.: ``The spectrum of a quasi-periodic Schr\"odinger operator'', {\it Commun. Math. Phys.} {\bf 111} (1987), 409-415.

\bibitem{TKN2} Thouless~D.~J., Kohmoto~M., Nightingale~M., den Nijs~M., ``Quantum Hall Conductance in Two Dimensional Periodic Potential'', {\it Phys. Rev. Lett.}, {\bf 49} (1982), 405-408.

\bibitem{Wil74} Williams R.F., ``Expanding attractors", {\it Publ. IHES}, {\bf 43} (1974), 169-203.

\end{thebibliography}
\end{document}